\documentclass[10pt]{article}

\usepackage[english,francais]{babel}

\usepackage{amsmath}
\usepackage{amssymb}

\usepackage{amsthm,amscd}

\usepackage{latexsym}

\usepackage{amsfonts}

\newtheorem{theoreme}{Th\'eor\`eme}[section]
\newtheorem{lemme}[theoreme]{Lemme}
\newtheorem{proposition}[theoreme]{Proposition}
\newtheorem{corollaire}[theoreme]{Corollaire}
\newtheorem{definition}[theoreme]{D\'efinition}
\newtheorem{remarque}[theoreme]{Remarque}

\theoremstyle{remark}
\numberwithin{equation}{section}

\author{Houmem Belkhechine\\
Facult\'e des Sciences de Gab\`es \\Tunisie\\
{\tt houmem@gmail.com}\and
Imed Boudabbous\\
Institut Pr\'eparatoire aux \'Etudes d'Ing\'enieurs de Sfax \\
Tunisie\\
{\tt imed.boudabbous@gmail.com}\and
Mohamed Baka Elayech \\
Institut Pr\'eparatoire aux \'Etudes d'Ing\'enieurs de Sfax \\
Tunisie\\
{\tt mohamedbaka.elayech@gmail.com}}
 \title{Les graphes (-1)-critiques}

\date{}

\begin{document}

\maketitle

\begin{abstract}
\selectlanguage{francais} \noindent \textit{\'Etant donn\'e un graphe (orient\'e) $G= ( S,A)$, une partie $X$ de $S$ est un intervalle de $G$
lorsque pour tous $a, b \in X$ et $x\in S-X$, $(a,x)\in A$ si et seulement si
$(b,x)\in A$ et $(x,a) \in A$ si et seulement si $(x,b)\in A$. Par exemple, $%
\emptyset$, $\{x\}$ $(x \in S)$ et $S$ sont des intervalles de $G$, appel\'es
intervalles triviaux. Un graphe, dont tous les intervalles sont triviaux,
est ind\'ecomposable; sinon, il est d\'ecomposable. Un sommet $x$ d'un
graphe ind\'ecomposable $G$ est critique si le graphe $G- x$ est
d\'ecomposable. En 1993, J.H. Schmerl et W.T. Trotter ont caract\'eris\'e
les graphes ind\'ecomposables dont tous les sommets sont critiques,
appel\'es graphes critiques. Dans cet article, nous caract\'erisons les
graphes ind\'ecomposables qui admettent un unique sommet non critique, que
nous appelons graphes (-1)-critiques}. Nous r\'epondons ainsi \`a une
question pos\'ee par Y. Boudabbous et P. Ille dans un article r\'ecent
\'etudiant les sommets critiques dans un graphe ind\'ecomposable.
\end{abstract}
\selectlanguage{english}
\begin{abstract}
  \noindent \textbf{The (-1)-critical graphs.
} \textit{Given a (directed) graph $G=(V,A)$, a subset $X$ of $V$ is
an interval of $G$ provided that for any $a, b\in X$ and $x\in V-X$,
$(a,x)\in A$ if and only if $(b,x)\in A$ and
$(x,a)\in A$ if and only if $(x,b)\in A$. For example, $\emptyset$, $%
\{x\}$ $(x \in V)$ and $V$ are intervals of $G$, called trivial intervals. A
graph, all the intervals of which are trivial, is indecomposable; otherwise,
it is decomposable. A vertex $x$ of an indecomposable graph is critical if $%
G-x$ is decomposable. In 1993, J.H. Schmerl and W.T. Trotter characterized
the indecomposable graphs, all the vertices of which are critical, called
critical graphs. In this article, we characterize the indecomposable graphs
which admit a single non critical vertex, that we call (-1)-critical graphs.}
This gives an answer to a question asked by Y. Boudabbous and P. Ille in a
recent article studying the critical vertices in an indecomposable graph.
\end{abstract}




\noindent {\bf Mots cl\'es:} intervalle, graphe
ind\'ecomposable, sommet critique, graphe (-1)-critique, graphe
d'ind\'ecomposabilit\'e.


\selectlanguage{french}

\section{Introduction}

\subsection{G\'en\'eralit\'es}

Un \textit{graphe} (\textit{orient\'e}) $G = (S(G), A(G))$ ou $(S, A)$, est
constitu\'e d'un ensemble fini $S$ de \textit{sommets} et d'un ensemble $A$
de couples de sommets distincts, appel\'es \textit{arcs} de $G$. L'\emph{%
ordre} (ou le \textit{cardinal}) du graphe $G$ est le nombre de ses sommets.
Soit $G= (S, A)$ un graphe.  \`A chaque partie $X$ de $S$ est
associ\'e le \textit{sous-graphe} $G(X) = (X, A \cap (X \times X))$ de $G$
\textit{induit} par $X$. Pour $X \subseteq S$ (resp. $x \in S$), le graphe $%
G(S-X)$, o\`u $S-X = \{s \in S : s \notin X\}$, (resp. $G(S-\{x\}$) est
not\'e $G-X$ (resp. $G-x$). Pour tous
sommets distincts $x$, $y$ de $S$, $x \longleftrightarrow y$
signifie $(x, y) \in A$ et $(y, x) \in A$; $x -- y$ signifie $%
(x, y) \notin A$ et $(y, x) \notin A$; $x \longrightarrow y$ signifie $(x,
y) \in A$ et $(y, x) \notin A$. Pour $x\in S$ et $Y\subseteq S$, $x
\longrightarrow Y$ signifie $x \longrightarrow y$ pour tout $y\in Y$.
Pour $X, Y\subseteq S$, $X \longrightarrow Y$ signifie $x \longrightarrow Y$
pour tout $x\in X$. D'une mani\`ere analogue, on d\'efinit pour $x\in S$ et
pour $X, Y\subseteq S$, $Y\longleftrightarrow x$, $x -- Y$, $X
\longleftrightarrow Y$ et $X -- Y$. Par exemple un \textit{tournoi} $T$ est
un graphe tel que pour tous $x \neq y$ dans $S(T)$, ou bien $x \longrightarrow y$
ou bien $y \longrightarrow x$. Un tournoi $T$ est un \textit{ordre total}
(ou une \textit{cha\^ine}), lorsque pour tous $x, y, z\in S(T) $, si $x
\longrightarrow y$ et $y \longrightarrow z$, alors $x \longrightarrow z$.
Dans un ordre total, la notation $x < y$ signifie $x \longrightarrow y$.
L'ordre total usuel $0 < \cdots < n$ est not\'e $O_{n}$.

Un graphe \textit{sym\'etrique} est un graphe $G$ tel que pour tous $x\neq
y$ dans $S(G)$, si $(x,y)\in A(G)$, alors $(y,x)\in A(G)$. \'Etant donn\'e un
sommet $x$ d'un graphe sym\'etrique $G$, un sommet $y\in S(G)$ est un
\textit{voisin} de $x$ (dans $G$), si $(x,y)\in A(G)$. On note $V_{G}(x)$
l'ensemble des voisins de $x$ dans $G$. le degr\'e de $x$ est $d_{G}(x)=
\mid V_{G}(x) \mid$. Lorsque $V_{G}(x)=\emptyset$, on dit que $x$ est un
sommet \textit{isol\'e} de $G$. Un graphe sym\'etrique $G$ est \textit{\
complet} (resp. \textit{vide}) lorsque pour tous $x\neq y$ dans  $S(G)$, $%
(x,y)\in A(G)$ (resp. $(x,y)\not\in A(G)$).

\'Etant donn\'e un graphe sym\'etrique $G$, une relation d'\'equivalence $%
\mathcal{R}$ est d\'efinie sur $S(G)$ comme suit. Pour tous $x \neq y$ dans
$S(G)$, $x$ $\mathcal{R}$ $y$ s'il existe une suite $x_{0} = x, \cdots,
x_{n}=y $ de sommets de $G$ telle que pour tout $i \in \{0, \cdots, n-1\}$, $%
(x_{i}, x_{i+1}) \in A(G)$. Les classes d'\'equivalence de $\mathcal{R}$
sont les \textit{composantes connexes} de $G$. Le graphe $G$ est \textit{%
connexe} lorsqu'il admet une seule composante connexe.

\'Etant donn\'es deux graphes $G = (S, A)$ et $G^{\prime}= (S, A^{\prime})$,
une bijection $f$ de $S$ sur $S^{\prime}$ est un \textit{isomorphisme} de $G$
sur $G^{\prime}$ si pour tous $x$, $y \in S$, $(x, y) \in A$ si et seulement
si $(f(x), f(y)) \in A^{\prime}$. Lorsqu'un tel isomorphisme existe, on dit
que $G$ et $G^{\prime}$ sont \textit{isomorphes}, et on note $G \simeq
G^{\prime}$. On dit qu'un graphe $G$ \textit{abrite} un graphe $H$ lorsque $%
H $ est isomorphe \`a un sous- graphe de $G$.

\`A chaque graphe $G$ est associ\'e son \textit{dual} $G^{\star}$ et son
\textit{compl\'ementaire} $\overline{G}$ d\'efinis sur $S(G)$ comme suit:
pour tous $x\neq y$ dans $S(G)$, $(x,y)\in A(G^{\star})$ si $(y,x)\in A(G)$; $%
(x,y)\in A(\overline{G})$ si $(x,y)\notin A(G)$. Un graphe $G$ est \textit{%
autodual} lorsqu'il est isomorphe \`a $G^{\star}$.

\subsection{Graphes ind\'ecomposables}

\'Etant donn\'e un graphe $G=(S, A)$, on introduit une relation
d'\'equivalen\-ce $\equiv_{G}$ (ou $\equiv$) sur l'ensemble des couples de
sommets distincts de $G$, d\'efinie comme suit: pour $x\neq y$ dans $S$ et $%
u \neq v$ dans $S$, $(x, y)\equiv_{G} (u,v)$ (ou $(x, y)\equiv (u,v)$ ) si $%
x\longrightarrow y$ et $u \longrightarrow v$, ou bien $y\longrightarrow x$ et $%
v\longrightarrow u$, ou bien $x \longleftrightarrow y$ et $u
\longleftrightarrow v$, ou bien $x -- y$ et $u -- v$. Dans le cas contraire
on note $(x, y)\not\equiv_{G} (u,v)$ (ou $(x, y)\not\equiv (u,v)$). Pour $%
x\in S$ et $Y\subseteq S - \{x\}$, $x \sim Y$ signifie que pour tous $y$, $z
\in Y$, $(x,y) \equiv (x,z)$. Pour $X, Y \subseteq S$, avec $X \cap Y =
\emptyset$, $X \sim Y$ signifie que pour tous $x$, $x^{\prime}\in X$ et pour
tous $y$, $y^{\prime}\in Y$, $(x,y) \equiv (x^{\prime}, y^{\prime})$. Par
ailleurs, une partie $I$ de $S$ est un \emph{intervalle} [5, 6, 8] (ou un
\textit{clan} \cite{E}) de $G$ lorsque pour tout $x\in S-I$, $x \sim I$. Par
exemple, $\emptyset $, $\{x\}$ o\`u $x \in S$, et $S$ sont les intervalles
\textit{triviaux} de $G$. Un graphe est \textit{ind\'ecomposable} [6, 8] (ou
\textit{primitif} \cite{E} ) si tous ses intervalles sont triviaux; il est
\textit{d\'ecomposable} dans le cas contraire. Nous introduisons quelques
notations et nous rappelons quelques propri\'et\'es des graphes
ind\'ecomposables.

\begin{definition}
Soit $G = (S, A)$ un graphe. \`A toute partie $X$ de $S$ telle que $\mid X
\mid \geq 3$ et $G(X)$ est ind\'ecomposable, on associe les
parties de $S-X$ suivantes:

\begin{itemize}
\item $[X] = \{x \in S-X: X$ est un intervalle de $G(X \cup \{x\})\}$.

\item Pour tout $u \in X$, $X(u) = \{x \in S-X: \ \{u, x\}$ est un
intervalle de $G(X \cup \{x\})\}$.

\item $Ext(X) = \{x \in S-X: \ G(X \cup \{x\})$ est ind\'ecomposable$\}$.
\end{itemize}

La classe form\'ee par $Ext(X)$, $[X]$ et $X(u)$, o\`u $u \in X$, est not\'ee $p_{X}$.
\end{definition}

\begin{lemme}
{\cite{E}} Soient $G = (S, A)$ un graphe et $X$ une partie de $S$ tels que $%
\mid X \mid \geq 3$ et $G(X)$ est ind\'ecomposable. La classe $p_{X} =
\{X(u): \ u \in X\} \cup \{Ext(X), [X]\}$ forme une partition de $S-X$. De
plus, les assertions suivantes sont v\'erifi\'ees.

\begin{itemize}
\item Soient $u \in X$, $x \in X(u)$ et $y \in S-(X \cup X(u))$. Si $G(X
\cup \{x, y\})$ est d\'ecomposable, alors $\{u, x\}$ est un intervalle de $%
G(X \cup \{x, y\})$.

\item Soient $x \in [X]$ et $y \in S-(X \cup [X])$. Si $G(X \cup \{x, y\})$
est d\'ecomposable, alors $X \cup \{y\}$ est un intervalle de $G(X \cup \{x,
y\})$.

\item Soient $x \neq y$ dans $Ext(X)$. Si $G(X \cup \{x, y\})$ est
d\'ecomposable, alors $\{x, y\}$ est un intervalle de $G(X \cup \{x, y\})$.
\end{itemize}
\end{lemme}

Le r\'esultat suivant est une cons\'equence directe du lemme 1.2.

\begin{corollaire}
{\cite{E}} Soit $G = (S, A)$ un graphe ind\'ecomposable. Si $X$ est une
partie de $S$ telle que $\mid X \mid \geq 3$, $\mid S-X \mid \geq 2$ et $G(X)
$ est ind\'ecomposable, alors il existe deux sommets distincts $x$ et $y$ de
$S-X$ tels que $G(X\cup\{x, y\})$ est ind\'ecomposable.
\end{corollaire}

Rappelons enfin le lemme suivant.

\begin{lemme}
{\cite{I}} Si $G = (S,A)$ est un graphe ind\'ecomposable avec $\mid S \mid \geq
5 $ et si $a \in S$, alors il existe une partie $X$ de $S$ telle que $\mid X
\mid = 4$ ou $5$, $a \in X$ et $G(X)$ est ind\'ecomposable.
\end{lemme}

\subsection{Graphes critiques}

Soit $G$ un graphe ind\'ecomposable. Un sommet $x$ de $G$ est
\textit{critique} si le graphe $G-x$ est d\'ecomposable. Lorsque tous les sommets de $G$ sont critiques, on dit que
$G$ est un \textit{graphe critique}. On g\'en\'eralise cette d\'efinition en disant
que le graphe $G$ est \textit{(-k)-critique} lorsqu'il admet exactement $k$
sommets non critiques. Notons que si un graphe $G$ est ind\'ecomposable,
alors $G^{\star}$ et $\overline{G}$ sont aussi ind\'ecomposables et ils ont
les m\^emes sommets critiques que $G$. Lorsqu'un sommet $a$ est l'unique sommet non
critique d'un graphe (-1)-critique $G$, on dit que le graphe $G$ est
(-1)-critique en $a$. J.H. Schmerl et W.T. Trotter \cite{ST} ont
caract\'eris\'e les graphes critiques. Nous rappelons cette
caract\'erisation dans le cas des tournois. Pour tout entier naturel $m$,
nous posons $N_{m} = \{0, \cdots, m\}$. Les tournois critiques sont, \`a des
isomorphismes pr\`es, les tournois $T_{2n+1}$, $U_{2n+1}$ et $V_{2n+1} $
d\'efinis sur $N_{2n}$, o\`u $n \geq 2$, comme suit.

\begin{itemize}
\item $T_{2n+1}(N_{n}) = 0 < \cdots < n$, $T_{2n+1}(\{n+1,\cdots,2n \})= n+1
< \cdots < 2n$, et pour tout $i \in N_{n-1}$, $\{i+1, \ldots,
n\}\longrightarrow i+n+1\longrightarrow N_{i}$.

\item $U_{2n+1}(N_{n}) = 0 < \cdots < n$, $(U_{2n+1})^{\star}(
\{n+1,\cdots,2n\})= n+1 < \cdots < 2n$, et pour tout $i \in N_{n-1}$, $%
\{i+1, \ldots, n\}\longrightarrow i+n+1 \longrightarrow N_{i}$.

\item $V_{2n+1}(N_{2n-1}) = 0 < \cdots < 2n-1$ et $\{2i+1: 0 \leq i \leq
n-1\} \longrightarrow 2n \longrightarrow \{2i: 0 \leq i \leq n-1\}$.
\end{itemize}

Dans cet article, nous caract\'erisons les graphes (-1)-critiques,
r\'epondant ainsi, \`a une question pos\'ee par Y. Boudabbous et P. Ille
\cite{BI}, et g\'en\'eralisant une r\'ecente caract\'erisation des tournois
(-1)-critiques \cite{BBD}.

\section{Caract\'erisation des graphes (-1)-critiques}

\subsection{Graphe d'ind\'ecomposabilit\'e}

La notion de \textit{graphe d'ind\'ecomposablit\'e} a \'et\'e introduite par
P. Ille [2, 7]. \`A chaque graphe $G = (S, A)$ est
associ\'e son graphe d'ind\'ecomposabilit\'e $I(G)$ d\'efini sur $S$ comme
suit. Pour tous $x \neq y$ dans $S$, $(x, y)$ est un arc de $I(G)$ si $G - \{x,
y\}$ est ind\'ecomposable. Notons que $I(G)$ est un graphe sym\'etrique et
que $I(\overline{G}) = I(G^{\star}) = I(G)$. Nous rappelons le lemme suivant.

\begin{lemme}
{\cite{BI}} Soient $G = (S,A)$ un graphe ind\'ecomposable d'ordre $\geq 5 $
et $x$ un sommet critique de $G$. Alors $\mid V_{I(G)}(x) \mid \leq 2$ et on
a:

\begin{itemize}
\item Si $V_{I(G)}(x) = \{y\}$, o\`u $y \in S$, alors $S - \{x, y\}$ est un
intervalle de $G-x$.

\item Si $V_{I(G)}(x) = \{y, z\}$, o\`u $y \neq z$ dans $S$, alors $\{y, z\}$
est un intervalle de $G-x$.
\end{itemize}
\end{lemme}

Le graphe d'ind\'ecomposabilit\'e est un outil important dans notre
construction des graphes (-1)-critiques. Afin de d\'ecrire les diff\'erents
graphes d'ind\'eco\-mposabilit\'e possibles d'un graphe (-1)-critique, nous
introduisons les graphes suivants. Le \textit{chemin} $P_{n}$ est le graphe
sym\'etrique d\'efini sur $N_{n}$ comme suit: Pour tous $i$, $j \in N_{n}$, $%
i \longleftrightarrow j$ si $\mid i- j\mid= 1$. Pour $n \geq 2 $, le \textit{%
cycle} $C_{n}$ est le graphe sym\'etrique obtenu \`a partir de $P_{n}$ en
ajoutant les arcs $(0, n)$ et $(n, 0)$. Tout graphe isomorphe \`a $P_{n}$
(resp. $C_{n}$) est appel\'e chemin (resp. cycle). La \textit{longueur} d'un
chemin, ou d'un cycle, est le nombre de paires $\{x, y\}$ de ses sommets
tels que $x \longleftrightarrow y$. Les \textit{extr\'emit\'es} (resp.
sommets \textit{\ internes}) d'un chemin sont ses sommets de degr\'e 1
(resp. de degr\'e 2). Un \textit{arbre} est un graphe sym\'etrique connexe
sans cycle. Les \textit{feuilles} d'un arbre sont ses sommets de degr\'e 1.
Un arbre \textit{\ \'etoil\'e} est un arbre $\mathcal{A}$ admettant un
unique sommet $a$ tel que $d_{\mathcal{A}}(a) \geq 3$, appel\'e \textit{%
source} de $\mathcal{A}$. Un arbre \textit{$a$-\'etoil\'e} est un arbre
\'etoil\'e de source $a$. \'Etant donn\'e un arbre \'etoil\'e $\mathcal{A}$,
une \textit{branche} de $\mathcal{A}$ est un chemin de $\mathcal{A}$ dont
les extr\'emit\'es sont la source et une feuille. Le degr\'e de $\mathcal{A}$
est le degr\'e de sa source, ou encore le nombre de ses branches ou de ses
feuilles. Nous consid\'erons enfin le graphe $R_{2n+1}$ d\'efini sur $N_{2n}$%
, o\`u $n \geq 2$, comme suit: $\{1, 3, \cdots ,2n-1\}\longrightarrow
2n\longrightarrow\{0, 2, \cdots, 2n -2\}$ et pour tous $x \neq y \in N_{2n-1}
$, $(x, y)$ est un arc de $R_{2n+1}$ si $x< y$ et si $x$ est impair ou $y$
est pair.

\begin{remarque}
{\cite{BI}} Le graphe $R_{2n+1}$, o\`u $n \geq 2$, est autodual et
(-1)-critique en $2n$. De plus, $I(R_{2n+1}) - 2n = P_{2n-1}$ et $2n$ est un
sommet isol\'e de $I(R_{2n+1})$.
\end{remarque}

\begin{lemme}
{\cite{BI}} Le graphe d'ind\'ecomposabilit\'e d'un graphe $G$ d'ordre $\geq
7 $ et (-1)-critique en $a$, admet une unique composante connexe de cardinal $\geq 2$. De plus, si $a$ est un sommet isol\'e de $I(G)$,
alors $G$ est isomorphe \`a $R_{2n+1}$ ou \`a $\overline{R_{2n+1}}$.
\end{lemme}

Nous compl\'etons ce lemme comme suit.

\begin{corollaire}
Soit $G$ un graphe d'ordre $\geq 7$, (-1)-critique en $a$ et tel que $a$
n'est pas un sommet isol\'e de $I(G)$. Soit $\mathcal{C}$ la composante
connexe non r\'eduite \`a un singleton de $I(G)$ et soit $X$ une partie de $%
S(G)$ telle que $G(X)$ est ind\'ecomposable. Si $\mathcal{C} \subseteq X$,
alors $S(G) = X$.
\end{corollaire}

\noindent\emph{Preuve\/. }
Par l'absurde, supposons $S(G) \neq X$. En appliquant plusieurs fois le corollaire 1.3 \`a partir de $G(X)$, on
obtient deux sommets, distincts ou non, $x,y \in S(G)- \mathcal{C}$ tels que
$G-\{x,y\}$ est ind\'ecomposable. Si $x=y$, $x$ est un sommet non critique
de $G$, contradiction car $a$ est l'unique sommet non critique de $G$. Si $x
\neq y$, $(x,y)$ est un arc de $I(G)$, contradiction car $x$ et $y$ sont des
sommets isol\'es de $I(G)$.
{\hspace*{\fill}$\Box$\medskip}

Ces r\'esultats nous am\`enent \`a associer \`a chaque graphe (-1)-critique $%
G$ d'ordre $\geq 7$, le sous-graphe $I^{\prime}(G)$ de $I(G)$, induit par sa
composante connexe non r\'eduite \`a un singleton.

\begin{proposition}
\'Etant donn\'e un graphe $G$ d'ordre $\geq 7$ et (-1)-critique en $a$,
l'une des assertions suivantes est v\'erifi\'ee:

\begin{itemize}
\item $I(G)$ est un cycle de longueur impaire.

\item $I^{\prime}(G)$ est un chemin de longueur $\geq 2$.

\item $I^{\prime}(G)$ est un arbre $a$-\'etoil\'e dont toutes les branches
sont de longueurs $\geq 2$ et admettant au plus une branche de longueur
impaire, et alors cette branche est de longueur $\geq$ 3.
\end{itemize}
\end{proposition}

\noindent\emph{Preuve\/. }
Soit $G$ un graphe d'ordre $\geq 7$ et (-1)-critique en $a$. Si $a$ est un
sommet isol\'e de $I(G)$ alors, d'apr\`es le lemme 2.3, $G$ est isomorphe
\`a $R_{2n+1}$ ou \`a $\overline{R_{2n+1}}$ et, d'apr\`es la remarque 2.2, $%
I^{\prime}(G)$ est un chemin de longueur $\geq 5$. Si $a$ n'est pas un sommet isol\'e de $I(G)$ alors $a$ est un sommet de $I^{\prime}(G)$. Supposons d'abord que $I(G)$ abrite un
cycle. Il existe alors une partie $X$ de $S(G)$ telle que $I(G)(X)$ est un
cycle. Posons $I(G)(X) = C_{m}$, o\`u $m \geq 2$. Le sommet $a$ appartient
\`a $X$, sinon il existe un sommet $\alpha \in X$ tel $d_{I(G)}(\alpha) \geq
3$. Comme $\alpha$ est un sommet critique de $G$, ceci contredit le lemme
2.1. On peut alors supposer que $a=0$. D'apr\`es le lemme 2.1, $\{0,
m-1\}$ est un intervalle de $G-m$. On a $m$ est pair, autrement, par le
lemme 2.1, $\{i, i+2\}$ est un intervalle de $G - \{i+1\}$ pour tout $i \in \{0, \cdots, m-3\}$, par suite $(m,0) \equiv (m,2) \equiv \cdots \equiv(m,m-1)$. Ainsi $%
\{0, m-1\}$ est un intervalle non trivial de $G$, contradiction. De plus, $%
S-X = \emptyset$, sinon un raisonnement analogue au pr\'ec\'edent utilisant le lemme 2.1 donne pour tout $\alpha \in S-X$, $(\alpha,1)
\equiv (\alpha,3) \equiv \cdots \equiv (\alpha,m-1) \equiv (\alpha,0) \equiv
(\alpha,2) \equiv \cdots \equiv (\alpha,m)$. Ainsi $X$ est un intervalle non
trivial de $G$, contradiction. Il s'ensuit que $I(G)$ est un cycle de
longueur impaire. Supposons \`a pr\'esent que $I(G)$ n'abrite pas un cycle,
c'est-\`a-dire que $I^{\prime}(G)$ est un arbre. Si $d_{I(G)}(a) \leq 2$
(resp. $d_{I(G)}(a) \geq 3$) alors, d'apr\`es le lemme 2.1, $I^{\prime}(G)$
est un chemin (resp. un arbre \'etoil\'e de source $a$). Supposons d'abord
que $I^{\prime}(G)$ est un chemin. D'apr\`es le lemme 1.4, il existe une
partie $X$ de $S(G)$ telle que $\mid X \mid = 4$ ou $5$, $a \in X$ et $G(X)$
est ind\'ecomposable. En appliquant plusieurs fois le corollaire 1.3 \`a
partir de $G(X)$, on obtient deux sommets $x,y \in S(G)$ tels que $G-\{x,y\}$
est ind\'ecomposable. Puisque $x$ est un sommet critique de $G$, alors $x
\neq y$. Il s'ensuit que $(x,y)$ est un arc de $I^{\prime}(G)$. Comme de
plus $a$ n'est pas un sommet isol\'e de $I(G)$, alors le chemin $%
I^{\prime}(G)$ est de longueur $\geq 2$. Nous posons maintenant, pour tous
entiers $h$, $l \geq 1$, $S_{h_{l}} = \{h_{0}, \cdots, h_{l}\}$, $h_{0} = a
= 0$, et nous notons par $P_{h_{l}}$ le chemin d\'efini sur $S_{h_{l}}$ par $%
A(P_{h_{l}}) = \{(h_{u}, h_{v}): \mid u - v \mid = 1\}$. Supposons par
l'absurde que $I^{\prime}(G)$ est un arbre \'etoil\'e admettant deux
branches distinctes $P_{1_{2p+1}}$ et $P_{2_{2q+1}}$ de longueurs impaires,
o\`u $p$, $q \in \mathbb{N}$. Une suite d'applications du lemme 2.1 donne $(0,1_{1}) \equiv (0, 1_{2p+1}) \equiv (1_{2p},1_{2p+1}) \not\equiv
(1_{2p},2_{1}) \equiv (0,2_{1})$. Soit $x \in S(G)- (S_{1_{2p+1}} \cup
S_{2_{2q+1}})$. Encore par le lemme 2.1, $(0,x) \equiv (1_{2p},x) \equiv
(1_{2p},2_{1}) \equiv (0,2_{1})$. Ainsi $(0,x) \equiv (0,2_{1})$. De
m\^eme, $(0,x) \equiv (0,1_{1})$. Contradiction car $(0,1_{1}) \not \equiv
(0,2_{1})$. Supposons que $I^{\prime}(G)$ est un arbre $0$-etoil\'e
admettant une branche $P_{1_{1}}$ de longueur 1. Consid\'erons deux autres
branches distinctes $P_{2_{2r}}$ et $P_{3_{2s}}$ de $I^{\prime}(G)$ o\`u $r$%
, $s \in \mathbb{N}^{\star}$. D'apr\`es le lemme 2.1, $0 \sim
(S-\{0,1_{1}\})$, et pour tous $x\in S-(\{1_{1}\}\cup S_{2_{2r}})$, $%
l\in\{1,\cdots,r\}$, on a $(2_{2l},x)\equiv(0,x)\equiv(0,3_{2s-1})
\equiv(2_{2r-1},3_{2s-1})\equiv(2_{2r-1},x)\equiv(2_{2l-1},x)$. Il s'ensuit
que $S_{2_{2r}}-\{0\}$ est un intervalle non trivial du graphe
ind\'ecomposable $G - \{0,1_{1}\}$, contradiction. D'o\`u, si elle existe,
la branche de longueur impaire de $I^{\prime}(G)$ est de longueur $\geq 3$.
{\hspace*{\fill}$\Box$\medskip}

La proposition 2.5 nous am\`ene \`a la distinction suivante des
graphes (-1)-critiques.

\subsection{Les graphes (-1)-critiques $G$ tels que $I(G)$ est un cycle}

Pour tout entier $p\geq 1$, nous consid\'erons le graphe $H_{2p+1}$ d\'efini
sur $N_{2p}$ comme suit: pour tous $x \neq y \in N_{2p}$
, $(x, y)$ est un arc de $H_{2p+1}$ si: ou bien $x< y$, $x$ est pair et $y$
est impair; ou bien $x> y$ et $x$ et $y$ sont de m\^eme parit\'e. Notons que
$H_{2p+1}$ est autodual en consid\'erant la permutation $\sigma$ d\'efinie par $\sigma(0) = 0$ et $\sigma(i) = 2p+1-i$ pour $i \neq 0$.

\begin{proposition}
\`A des isomorphismes pr\`es, les graphes (-1)-critiques d'ordre $\geq 7$ et
dont le graphe d'ind\'ecomposabilit\'e est un cycle sont $H_{2p+1}$ et $%
\overline{H_{2p+1}}$, o\`u $p \geq 3$. De plus, $0$ est le sommet non
critique de $H_{2p+1}$.
\end{proposition}

\noindent\emph{Preuve\/. } Nous commen\c{c}ons par \'etablir que pour tout $%
p \geq 2$, $H_{2p+1}$ est (-1)-critique en $0$ et que $I(H_{2p+1})=C_{2p+1} $%
. Montrons d'abord, par r\'ecurrence, que pour tout $p \geq 1$, $H_{2p+1}$
est ind\'ecomposable. Il est clair que $H_{3}$ est ind\'ecomposable. Soit
maintenant $p\geq 2$. Le graphe $H_{2p-1}$ \'etant ind\'ecompo\-sable par
hypoth\`ese de r\'ecurrence, on v\'erifie que $H_{2p+1}$ est
ind\'ecomposable en consid\'erant la partition $p_{W}$, o\`u $W = S(H_{2p-1})
$. Il suffit de constater que $2p\in W(2p-2)$, $2p-1\not\in W(2p-2)$ et $%
(2p-1, 2p)\not\equiv (2p-1, 2p-2)$. Tout sommet $i \neq 0$ de $H_{2p+1}$ est
critique car $\{i-1, i+1\}$ est un intervalle non trivial de $H_{2p+1}-i$.
En remarquant que $H_{2p+1} - \{0,1\} \simeq H_{2p-1}$, nous v\'erifions
maintenant que $H_{2p+1} - 0$ est ind\'ecomposable en consid\'erant la
partition $p_{Y}$, o\`u $Y = \{2, \cdots, 2p\}$. D'une part $1 \not\in [Y]$
car $2--1$ et $3 \longrightarrow 1$. D'autre part, pour tout $u \in Y$, $1
\not\in Y(u)$. En effet, $(2p-1, 1) \not \equiv (2p-1, u)$ si $u$ est pair; $%
2 \longrightarrow u$ et $2 -- 1$ si $u$ est impair. Il s'ensuit, d'apr\`es le lemme 1.2, que $1 \in
Ext(Y)$, c'est-\`a-dire $0$ est un sommet non critique de $H_{2p+1}$. Pour
v\'erifier que $I(H_{2p+1})=C_{2p+1}$, il suffit, d'apr\`es la proposition
2.5, de v\'erifier que pour tout $i \in N_{2p}$, $%
H_{2p+1} - \{i, i+1\}$ est ind\'ecomposable. Pour tout $i \in N_{2p-1}$, $%
H_{2p+1}-\{i, i+1\} \simeq H_{2p-1}$. De plus, $%
H_{2p+1}-\{2p, 0\} \simeq H_{2p-1}$. Comme $H_{2p-1} $ est ind\'ecomposable, il
s'ensuit que $I(H_{2p+1})=C_{2p+1}$.

Soit maintenant un graphe (-1)-critique $G$ d'ordre $\geq 7$ et dont le
graphe d'ind\'ecomposabilit\'e est un cycle. D'apr\`es la proposition 2.5,
le graphe $G$ est de cardinal impair. On pose alors $S(G) = N_{2p}$, o\`u $p \geq 3$, $0$ est le sommet non critique de $G$
et $I(G)=C_{2p}$. D'apr\'es le lemme 2.1, pour tout $i\in
\{1,\cdots,2p\}$ , $\{i-1, i+1\}$ est un intervalle de $G-i$, alors $(0,
1)\equiv (0, 3)\equiv\cdots\equiv (0, 2p-1)\equiv (2,
2p-1)\equiv\cdots\equiv (2p-2,2p-1)\equiv(2p-2, 0)\equiv(2p, 0) $. Il
s'ensuit que $(0, 1)\not\equiv(0,2p)$, autrement $(0, 1)\equiv (0,
3)\equiv\cdots\equiv (0, 2p-1)\equiv (0,2p) \equiv (0,2p-2) \equiv \cdots
\equiv (0,2)$, c'est-\`a-dire $\{1, \cdots, 2p\}$ est un intervalle non
trivial de $G$, contradiction. Ainsi $(0, 1)\not\equiv(1,0)$ et quitte \`a
remplacer $G$ par $G^{\star}$, on peut supposer que $0\longrightarrow 1$ et
donc $2p\longrightarrow 0$. Pour $i \leq j$ dans $N_{p-1}$, $2i \longrightarrow
2j+1$ car $0 \longrightarrow 1$ et $(2i,2j+1) \equiv \cdots\equiv(0, 2j+1)
\equiv \cdots\equiv(0, 1)$. Pour $i < j$ dans $N_{p}$, $2j \longrightarrow 2i$
car $2p \longrightarrow 0$ et $(2j, 2i) \equiv \cdots \equiv(2p, 2i ) \equiv
\cdots \equiv(2p, 0)$. Pour $i < j$ dans $N_{p-1} $, $2j+1\longrightarrow 2i+1$
car $0 \longrightarrow 1$ et $(2j+1, 2i+1) \equiv \cdots \equiv (2p-1, 2i+1
) \equiv \cdots \equiv(2p-1, 1) \equiv(0, 1) $. Pour $i < j$ dans $N_{p}$, $%
(2i+1, 2j) \equiv \cdots \equiv (1, 2j ) \equiv \cdots \equiv(1, 2)$, en particulier $(2i+1, 2p) \equiv (1,2)$. Or $%
(1, 2)\equiv (2,1)$ sinon, comme $0 \longrightarrow 1$ et $\{0,2\}$ est un
intervalle de $G-1$ et n'est pas un intervalle de $G$, alors $%
1\longrightarrow 2$, il s'ensuit que pour tout $i \in \{1,\cdots, p-1\}$, $%
2i+1\longrightarrow\{1, 2p\} \longrightarrow 2i$, en particulier $\{1, 2p\}$
est un intervalle non trivial de $G-0$. Contradiction. Ainsi, ou bien $1 -- 2
$ et dans ce cas $G= H_{2p+1}$, ou bien $1\longleftrightarrow 2$ et dans ce
cas $G =(\overline{ H_{2p+1}})^{\star} \simeq \overline{H_{2p+1}}$.
{\hspace*{\fill}$\Box$\medskip}

\subsection{Les graphes (-1)-critiques $G$ tels que $I^{\prime}(G)$ est un
chemin}

Nous construisons les graphes (-1)-critiques $G$ de ce paragraphe sur $S(G)
= N_{m}$, $N_{m} \cup \{\alpha\}$ ou $N_{m} \cup \{\alpha, \beta\}$, o\`u $%
\{\alpha, \beta\}$ est une paire d'\'el\'ements distincts tels que $%
\{\alpha, \beta\} \cap \mathbb{N} = \emptyset$.

Le corollaire suivant d\'ecoule directement de la remarque 2.2 et du lemme
2.3.

\begin{corollaire}
\`A des isomorphismes pr\`es, les graphes $G$ d'ordre $\geq 7$,
(-1)-critiques en $a$ et tels que $I^{\prime}(G)$ est un chemin et $a$ est
un sommet isol\'e de $I(G)$, sont $R_{2n+1}$ et $\overline{R_{2n+1}}$, o\`u $%
n \geq 3$.
\end{corollaire}

Nous distinguons maintenant les cas o\`u le sommet non critique est une
extr\'emit\'e ou un sommet interne de $I^{\prime}(G)$.

\subsubsection{Les graphes (-1)-critiques $G$ en une extr\'emit\'e de $%
I^{\prime}(G)$}

Nous introduisons, pour $m \geq 2$, la classe $\mathcal{F}_{m}$ des graphes $%
G$ d\'efinis sur $N_{m}$ et tels que $N_{m}-\{0, 1\}$ est un intervalle de $%
G-0$ et pour tout $i\in\{1,\cdots, m-1\}$, $\{i-1, i+1\}$ est un intervalle
de $G-i$ sans \^etre un intervalle de $G$. Notons d'abord la remarque
suivante.

\begin{remarque}
Soit $G \in \mathcal{F}_{m}$, o\`u $m \geq 2$. Pour tout $i \in
N_{m-1}$, $G -\{i, i+1\} \simeq G- \{m-1, m\}$. De plus, si $m \geq 3$, $G-m
\in \mathcal{F}_{m-1}$.
\end{remarque}

Nous caract\'erisons la classe $\mathcal{F}_{m}$ comme suit.

\begin{lemme}
\'Etant donn\'e un graphe $G$ d\'efini sur $N_{m}$, o\`u $m\geq 2 $. Alors $G\in
\mathcal{F}_{m}$ si et seulement si $(0,1) \not \equiv (2,1)$ et pour tous $x < y$ dans $N_{m}$, $(x,y) \equiv (1,2)$ si $x$ est impair, $(x,y) \equiv (0,2)
$ si $x$ et $y$ sont pairs, et $(x,y) \equiv (0,1)$ si $x$ est pair et $y$
est impair.
\end{lemme}

\noindent\emph{Preuve\/. }
Soit $G\in \mathcal{F}_{m}$, o\`u $m\geq 2$. On a $(0,1) \not \equiv (2,1)$ car $%
\{0,2\}$ est un intervalle de $G -1$ et n'est pas un intervalle de $G$. Soit
$x < y$ dans $N_{m}$. Pour $y \neq 1$, $(1, y)\equiv (1,2)$ car $N_{m}-\{0, 1\}$
est un intervalle de $G-0$. Comme de plus pour tout $i\in\{1,\cdots, m-1\}$,
$\{i-1, i+1\}$ est un intervalle de $G-i$, alors $(1, 2)\equiv (1,
y)\equiv\cdots\equiv (x, y)$ si $x$ est impair; $(0, 1) \equiv\cdots\equiv
(0,y) \equiv\cdots\equiv (x, y)$ si $x$ est pair et $y$ est impair, et $(0,
2) \equiv\cdots\equiv (0,y) \equiv\cdots\equiv (x, y)$ si $x$ et $y$ sont
pairs.

R\'eciproquement, soit $i\in \{1,\cdots, m-1\}$. Remarquons que pour $x\in
N_{m}-\{i-1, i, i+1\}$, $(x, i-1)\equiv (x, i+1)$, de sorte que $\{i-1, i+1\}
$ est un intervalle de $G-i$. Si $i $ est pair (resp. impair), alors $%
(i-1, i)\equiv ( 1, 2)$ et $(i+1, i)\equiv ( 1, 0)$ (resp. $(i-1, i)\equiv (
0, 1)$ et $(i+1, i)\equiv ( 2, 1)$). Puisque $(0,1) \not \equiv (2,1)$, $(i,
i-1) \not \equiv (i, i+1)$ et par suite $\{i-1, i+1\}$ n'est pas un
intervalle de $G$. Enfin, $N_{m}-\{0, 1\}$ est un intervalle de $G-0$ car
pour tout $y\in N_{m}-\{0, 1\}$, $(1, y)\equiv ( 1, 2)$.
{\hspace*{\fill}$\Box$\medskip}

Le r\'esultat suivant caract\'erise les graphes
ind\'ecomposables de $\mathcal{F}_{m}$.

\begin{lemme}
Soit $G \in \mathcal{F}_{m}$, o\`u $m \geq 2$. Alors $G$ est
ind\'ecomposable si et seulement si ou bien $m$ est pair et $(0,1)
\not\equiv (0,2) \not\equiv (1,2)$, ou bien $m$ est impair et $(0,2)
\not\equiv (0,1) \not\equiv (1,2)$.
\end{lemme}

\noindent\emph{Preuve\/. } Soit $G\in \mathcal{F}_{m}$. On suppose
que $(0,1) \not\equiv (0,2) \not\equiv (1,2)$ lorsque $m$ est pair, et que $%
(0,2) \not\equiv (0,1) \not\equiv (1,2)$ lorsque $m$ est impair. Montrons,
par r\'ecurrence sur $m$, que $G$ est ind\'ecomposale. En utilisant le lemme
2.9, on v\'erifie que $G$ est ind\'ecomposable pour $m=2$ et pour $m=3$.
Soit $m \geq 4$. D'apr\`es la remarque $2.8$, $G(N_{m-2}) \in \mathcal{F}%
_{m-2}$. Il s'ensuit que $G(N_{m-2})$ est ind\'ecomposable par hypoth\`ese
de r\'ecurrence. Comme $m \in N_{m-2}(m-2)$, $m-1 \notin N_{m-2}(m-2)$ et $%
(m-1, m-2) \not\equiv (m-1,m)$, alors $G$ est ind\'ecomposable d'apr\`es le
lemme 1.2. R\'eciproquement, si $m$ est pair et si $(0,2)\equiv (0,1)$ ou $%
(0,2)\equiv (1,2)$ (resp. si $m$ est impair et si $(0,1)\equiv (0,2)$ ou $%
(0,1)\equiv (1,2)$), on v\'erifie que $\{1, 2,\cdots, m\}$ ou $N_{m-1}$ est
un intervalle non trivial de $G$.
{\hspace*{\fill}$\Box$\medskip}

Afin de caract\'eriser les graphes (-1)-critiques de ce paragraphe, nous
introduisons la classe $\mathcal{F}$ des graphes $G = (S,A)$ tels que $S =
N_{m}$, $N_{m} \cup \{\alpha\}$ ou $N_{m} \cup \{\alpha, \beta\}$, o\`u $m
\geq 2$ $; G(N_{m}) \in \mathcal{F}_{m}$; $(0,1) \equiv (0,2)$ si et
seulement si $S - N_{m} \neq \emptyset$; pour tous $i \in N_{m}$ et $\gamma
\in S- N_{m}$, $(i, \gamma) \equiv (1,2)$ si $i$ est impair et $(i, \gamma)
\equiv (0, \gamma)$ si $i$ est pair; et tels que:

\begin{itemize}
\item Si $S = N_{m}$, alors $(0,2) \not\equiv (1,2) \not\equiv (0,1)$.

\item Si $S - N_{m} = \{\alpha\}$, alors $(0, \alpha) \equiv (\alpha, 0)$ si
$(0,1) \equiv (1,2)$; $(0, \alpha) \not\equiv (0,1)$ et $(0, \alpha)
\not\equiv (1,2)$ si $(0,1) \not\equiv (1,2)$.

\item Si $S - N_{m} = \{\alpha, \beta\}$, alors $(\beta,\alpha) \not\equiv
(1,2) \not\equiv (0,1)$, $(0, \alpha) \equiv (1,2)$ et $(0, \beta) \equiv
(0,1)$.
\end{itemize}

\begin{lemme}
\'Etant donn\'e un graphe $G$ d'ordre $\geq 4$ de la classe $\mathcal{F}$. Alors $G$ est (-1)-critique en $m$. De plus, si $G$ est d'ordre $\geq 7$, alors $I^{\prime}(G) = P_{m}$.
\end{lemme}

\noindent\emph{Preuve\/. } Soit $G = (S,A)$ un graphe de la classe $\mathcal{%
\ F}$. Supposons d'abord que $S = N_{m}$. Dans ce cas, $m \geq 3$ et,
d'apr\`es la remarque 2.8, $G(N_{m-1}) \in \mathcal{F}_{m-1}$. En utilisant
le lemme 2.10, les graphes $G$ et $G-m$ sont ind\'ecomposables. Par
d\'efinition de $\mathcal{F}_{m}$, pour tout $i \in N_{m-1}$, $i$ est un
sommet critique de $G$. Il s'ensuit que $G$ est (-1)-critique en $m$.
Lorsque $\mid S \mid\geq 7 $, le graphe $G-\{m-1, m\}$ est ind\'ecomposable
et pour tout $i \in N_{m-2}$, $i$ est un sommet critique de $G-m$. Il
s'ensuit que $V_{I(G)}(m) = \{m-1\}$ et que, par la proposition 2.5, $%
I^{\prime}(G)$ est un chemin. De plus, par la remarque 2.8, pour tout $i \in
N_{m-1}$, $G - \{i, i+1\} \simeq G - \{m-1,m\}$, en particulier le graphe $%
G - \{i, i+1\}$ est ind\'ecomposable. On conclut que $I^{\prime}(G) = I(G) =
P_{m}$. Supposons maintenant que $S - N_{m} \neq \emptyset$. On v\'erifie
que $G(S - \{2, \cdots, m\})$ est ind\'ecomposable. Soit $i \in \{2, \cdots,
m\}$. On pose $S_{i} = S- \{i, \cdots, m\}$. Montrons que si $%
G(S_{i})$ est ind\'ecomposable, il en est de m\^eme pour $G(S_{i}
\cup \{i\})$. On a $(i,1) \equiv (2,1)$ et $(i,0) \equiv (1,0)$. Si $i$ est
impair, alors, pour $\gamma \in S - N_{m}$, $(i,\gamma) \equiv (1,2) \not
\equiv (i,0)$. Si $i$ est pair et $(i,0) \equiv (i,1)$, alors $S - N_{m} =
\{\alpha\}$ et $(\alpha,0) \equiv (0,\alpha) \equiv (i,\alpha) \not \equiv
(i,0)$. Ainsi $i \not \in [S_{i}]$. De plus, $i \not \in S_{i}(0)$ car $%
(1,0) \not \equiv (1,i)$, et pour $j \in \{1, \cdots, i-1\}$, $i \not \in
S_{i}(j)$ sinon $G ((S_{i} - \{j\}) \cup \{i\})$ est ind\'ecomposable.
Contradiction car $\{j-1, j+1\}$ est un intervalle de $G - j$. Si $S - N_{m}
= \{\alpha\}$, $(i,0) \not \equiv (\alpha,0)$ et donc $i \not \in
S_{i}(\alpha)$. Si $S - N_{m} = \{\alpha, \beta\}$, alors $i \not \in
S_{i}(\alpha)$ car $(0,\alpha) \equiv (1,2) \not \equiv (0,1) \equiv (0,i)$.
De plus, $(\beta,\alpha) \not \equiv (1,2) \equiv (0,\alpha)\equiv (i,\alpha)
$, donc $i \not \in S_{i}(\beta)$. Nous concluons, par le lemme 1.2, que $i
\in Ext(S_{i})$. Il s'ensuit que les graphes $G$ et $G-m$ sont
ind\'ecomposables. Pour tout $i \in \{1, \cdots, m-1\}$, $i$ est un sommet
critique de $G$ car $\{i-1,i+1\} $ est un intervalle non trivial de $G-i$.
De m\^eme $0$ et $\alpha$ sont des sommets critiques de $G$ car $S - \{0,1\}$
et $S - \{0, \alpha\}$ sont des intervalles non triviaux respectifs de $G-0$
et de $G- \alpha$. Dans le cas o\`u $\mid S- N_{m} \mid = 2$, $N_{m}$ est un
intervalle non trivial de $G - \beta$. Le graphe $G$ est alors (-1)-critique
en $m$. De m\^eme, $G - \{m-1, m\}$ est ind\'ecomposable et les sommets de $%
G-\{m-1,m\}$ sont des sommets critiques de $G-m$ lorsque $m \geq 3$. Il
s'ensuit, d'apr\`es la proposition 2.5, que $I^{\prime}(G)$ est un chemin
lorsque $\mid S \mid \geq 7$. Remarquons que pour tout $i\in N_{m-2}$, $%
G-\{i, i+1\} \simeq G-\{m-1, m\}$, en particulier $G-\{i, i+1\}$ est
ind\'ecomposable. De plus, $\alpha$ est un sommet isol\'e de $I(G)$. Il en
est de m\^eme pour le sommet $\beta$ lorsque $\mid S - N_{m} \mid = 2$. On
conclut que $I^{\prime}(G) = P_{m}$ lorsque $\mid S \mid \geq 7 $.
{\hspace*{\fill}$\Box$\medskip}

\begin{proposition}
\`A des isomorphismes pr\`es, les graphes $G$ d'ordre $\geq 7$,
(-1)-critiques en $m$ et tels que $I^{\prime}(G) = P_{m}$, sont les graphes
d'ordre $\geq 7$ de la classe $\mathcal{F}$.
\end{proposition}

\noindent\emph{Preuve\/. } Soit $G = (S,A)$ un graphe d'ordre $\geq 7$,
(-1)-critique en $m$ et tel que $I^{\prime}(G)= P_{m}$. D'apr\`es la
proposition 2.5, $m \geq 2$. D'apr\`es le lemme 2.1, pour tout $i \in \{1,
\cdots, m-1\}$, $\{i-1, i+1\} $ est un intervalle de $G- i$ et donc un
intervalle de $G(N_{m}) - i$. De plus, $\{i-1, i+1\} $ n'est pas un
intervalle de $G$. Il s'ensuit que $(i, i+ 1)\not\equiv (i, i- 1)$, en
particulier, $\{i-1, i+1\}$ n'est pas un intervalle de $G(N_{m})$. D'apr\`es
le lemme 2.1, $S-\{0, 1\}$ est un intervalle de $G-0$, en particulier, $%
N_{m}-\{0, 1\}$ est un intervalle de $G(N_{m})-0$. On conclut que $%
G(N_{m})\in \mathcal{F}_{m}$. Si $S-N_{m} \neq\emptyset$, alors, encore par
le lemme 2.1, pour tous $\gamma \in S-N_{m} $ et $i\in N_{m}$, $(i,
\gamma)\equiv (0, \gamma)$ si $i$ est pair et $(i, \gamma)\equiv (1, 2)$ si $%
i$ est impair.

Nous montrons maintenant que $(0,1) \equiv (0,2)$ si et seulement si $S -
N_{m} \neq \emptyset$. Si $(0,1) \equiv (0,2)$, alors $S - N_{m} \neq
\emptyset$ sinon $\{1, \cdots, m\}$ est un intervalle non trivial de $G$,
contradiction. Supposons que $(0,1) \not\equiv (0,2)$. Il suffit de montrer
que $(1, 2) \not\equiv (0, 1)$ si $m$ est impair et $(1, 2) \not\equiv (0, 2)
$ si $m$ est pair. En effet, dans ce cas, $G(N_{m})$ est ind\'ecomposable
d'apr\`es le lemme 2.10. Il s'ensuit que $S = N_{m}$ d'apr\`es le corollaire
2.4. Supposons \`a pr\'esent que $(0,1) \equiv (1,2)$ et $m$ est impair.
Encore d'apr\`es le lemme 2.10 et la remarque 2.8, $G(N_{m})$ est
d\'ecomposable et $G(N_{m-1})$ est ind\'ecomposable. Il existe alors $\mu
\in S-N_{m}$ tel que $(0, \mu)\not\equiv (1, 2)$. En utilisant le lemme 1.2,
$G(N_{m}\cup\{\mu\})$ est ind\'ecomposable car $m \in [N_{m-1}]$, $%
\mu\not\in [N_{m-1}]$ et $(m, \mu)\not\equiv (m, m-1)$. D'apr\`es le
corollaire 2.4, $G = G(N_{m}\cup\{\mu\})$. Ainsi, $G-\{\mu, m\}$ est
ind\'ecomposable, contradiction car $\mu$ est un sommet isol\'e de $I(G)$.
Supposons enfin que $(0,2) \equiv (1,2)$ et $m$ est pair. D'apr\`es le lemme
2.10, $G(N_{m})$ est d\'ecomposable. Ainsi $S - N_{m} \neq \emptyset$.  De
plus, d'apr\`es le lemme 2.10 et la remarque 2.8,  $G(N_{m-1})$ est
ind\'ecomposable pour $m\geq 4$. Supposons qu'il existe un sommet $\gamma
\in S-N_{m}$ tel que $(2,1)\not \equiv (0, \gamma) \not \equiv (1,2)$. Comme
$G(\{0,1,2,\gamma\})$ est ind\'ecomposable alors, d'apr\`es le corollaire
2.4, $m\geq 4$. On obtient $G(N_{m} \cup \{\gamma\})$ est ind\'ecomposable
car $m \in [N_{m-1}]$, $\gamma \not \in [N_{m-1}]$ et $(m,\gamma) \not
\equiv (m,m-1)$. D'apr\`es le corollaire 2.4, $(\gamma, m)$ est un arc de $%
I(G)$, contradiction. Il s'ensuit que $(1,2) \not \equiv (2,1)$ et que $%
S-N_{m} = E_{1} \cup E_{2}$, o\`u $E_{1} = \{x \in S- N_{m}: (0, x)\equiv
(1, 2)\}$ et $E_{2} = \{x \in S- N_{m}: (0, x)\equiv (2, 1)\}$. Notons que $%
E_{1} \neq \emptyset$ (resp. $E_{2} \neq \emptyset$), sinon $S- \{m\}$
(resp. $N_{m}$) est un intervalle non trivial de $G$. Comme $E_{2} \cup N_{m}
$ n'est pas un intervalle de $G$, il existe $(e_{1}, e_{2}) \in E_{1} \times
E_{2}$, tel que $(e_{2}, e_{1}) \not \equiv (1,2)$. On v\'erifie que $%
G(\{0,1, e_{2}\})$ est ind\'ecomposable. En utilisant le lemme 1.2, $%
G(\{0,1,2,e_{1}, e_{2}\})$ est ind\'ecomposable car $2\in[\{0,1, e_{2}\}]$, $%
e_{1}\not\in [\{0,1, e_{2}\}]$ et $(e_{1},2)\not\equiv (e_{2},2)$. Ainsi,
par le corollaire 2.4, on obtient $m\geq 4$. Posons $X = N_{m-1} \cup
\{e_{1}, e_{2}\}$, nous montrons que $G(X)$ est ind\'ecomposable. On
v\'erifie que $e_{2} \in Ext(N_{m-1})$. En effet, $e_{2} \not \in [N_{m-1}]$
car $(0, e_{2}) \equiv (2,1) \not \equiv (1,2) \equiv (1, e_{2})$. De plus,
pour tout $i \in N_{\frac{m}{2}-1}$, $e_{2} \not \in N_{m-1}(2i)$ car $%
(e_{2}, m-1) \equiv (2,1) \not \equiv (0,1) \equiv (2i, m-1)$, et $e_{2}
\not \in X(2i+1)$ car $(e_{2}, 0) \equiv (1,2) \not \equiv (1,0) \equiv
(2i+1, 0)$. Comme de plus, $e_{1} \in [N_{m-1}]$ et $(e_{2}, e_{1}) \not
\equiv (1,2) \equiv (1, e_{1})$, alors $G(X)$ est ind\'ecomposable. On a $m
\in Ext(X)$. En effet, d'une part $m \not \in [X]$ car $(1,m) \not \equiv
(e_{1}, m)$. D'autre part, pour tout $i \in N_ {\frac{m}{2}-1}$, $m \not \in
X(2i)$ car $(2i, m-1) \equiv (0,1) \not \equiv (2,1) \equiv (m, m-1)$, et $m
\not \in X(2i+1)$ car $(2i+1, e_{2}) \equiv (1,2) \not \equiv (2,1) \equiv
(m, e_{2})$. De plus $m\not\in X(e_{1})\cup X(e_{2})$ car $%
(e_{2},m)\equiv(m,e_{1})\equiv(1,2)\not\equiv(e_{2},e_{1})$. Ainsi $G(N_{m}
\cup \{e_{1}, e_{2}\})$ est ind\'ecomposable et, d'apr\`es le corollaire
2.4, $G = G(N_{m} \cup \{e_{1}, e_{2}\})$. Il s'ensuit que $G- \{m, e_{1}\}$
est ind\'ecomposable, contradiction car $e_{1}$ est un sommet isol\'e de $%
I(G)$.

Si $S = N_{m}$, alors $G \in \mathcal{F}_{m}$ et $G-m \in \mathcal{F}_{m-1}$
et, d'apr\`es le lemme 2.10, $(0,2) \not\equiv (1,2) \not\equiv (0,1)$. Il
s'ensuit que $G \in \mathcal{F}$. Supposons alors que $S - N_{m} \neq
\emptyset$ . Dans ce cas $(0, 1)\equiv (0, 2)$. Supposons d'abord que $(0,
1)\not\equiv (1, 2)$. S'il existe $\omega \in S-N_{m}$ tel que $(0,
\omega)\not\equiv (0, 1)$ et $(0, \omega)\not\equiv (1, 2)$ alors, $G(N_{m}
\cup \{\omega\})$ est isomorphe \`a un graphe de la classe $\mathcal{F}$. Il
s'ensuit que $G(N_{m} \cup \{\omega\})$ est (-1)-critique et que, d'apr\`es
le corollaire 2.4, $G = G(N_{m} \cup \{\omega\})$. Sinon, $S -N_{m} = E_{1}
\cup E_{3}$ o\`u $E_{3} = \{x \in S- N_{m}: (0, x)\equiv (0, 1)\}$.
Remarquons que $E_{1} \neq\emptyset$ (resp. $E_{3} \neq\emptyset$), sinon $%
(E_{3} \cup N_{m})-\{0\}$ (resp. $N_{m}$) est un intervalle non trivial de $G
$, contradiction. Comme $E_{3} \cup N_{m} $ n'est pas un intervalle de $G$,
il existe $(e_{1}, e_{3}) \in E_{1} \times E_{3}$, tel que $(e_{3},
e_{1})\not\equiv (1, 2)$. Ainsi, $G(N_{m} \cup \{e_{1}, e_{3}\})$ est
isomorphe \`a un graphe de $\mathcal{F}$. Il s'ensuit que $G(N_{m}\cup
\{e_{1}, e_{3}\})$ est (-1)-critique et que, d'apr\`es le corollaire 2.4, $G
= G(N_{m} \cup \{e_{1}, e_{3}\})$. Supposons enfin que $(0, 1) \equiv (1, 2)$%
. D'apr\`es le lemme 2.9, d'une part $(0,1)\not \equiv (1,0)$ et d'autre
part, pour tous $x < y \in N_{m}$, $(x,y) \equiv (0,1)$. Il s'ensuit que $%
G(N_{m}) = O_{m}$ ou $O^{\star}_{m}$. Remarquons que si $G \in \mathcal{F}$,
alors $G^{\star} \in \mathcal{F}$. On peut alors supposer que $G(N_{m})=
O_{m}$. Il suffit de montrer qu'il existe un sommet $\nu \in S-N_{m}$, tel
que $(0, \nu)\not\equiv (1, 2)$ et $(0, \nu)\not\equiv (2, 1)$. En effet,
dans ce cas, $(0, \nu) \equiv (\nu, 0)$ et par suite, $G(N_{m} \cup \{\nu\})$
est isomorphe \` a un graphe de la classe $\mathcal{F}$. Il s'ensuit que $%
G(N_{m} \cup \{\nu\})$ est (-1)-critique et que, d'apr\`es le corollaire
2.4, $G = G(N_{m} \cup \{\nu\})$. Suppposons alors, par l'absurde, que $S -
N_{m} = E_{1} \cup E_{2}$. Remarquons que $E_{1}\neq\emptyset$ sinon, $E_{2}
\cup N_{m-2}$ est un intervalle non trivial de $G-m$ si $m$ est impair, $%
E_{2} \cup N_{m-1}$ est un intervalle non trivial de $G$ si $m$ est pair.
Contradiction car $G$ et $G - m$ sont ind\'ecomposables. De m\^eme $%
E_{2}\neq\emptyset$, sinon $N_{m}$ est un intervalle non trivial de $G$,
contradiction. L'entier $m$ est pair, autrement, pour $y \in E_{2}$, $%
G(N_{m} \cup \{y\}) \simeq V_{m+2}$. D'apr\`es le corollaire 2.4, $G =
G(N_{m} \cup \{y\}) \simeq V_{m+2}$, contradiction car $V_{m+2}$ est un
tournoi critique. Comme $N_{m} \longrightarrow E_{1}$ et $E_{2}\cup N_{m} $
n'est pas un intervalle de $G$, il existe $(b_{1},b_{2}) \in E_{1} \times
E_{2}$, tel que $(b_{2}, b_{1})\not\equiv (1, 2)$. En remarquant que $%
G(N_{m-1} \cup \{b_{2}\}) \simeq V_{m+1}$, on v\'erifie que $G(N_{m} \cup
\{b_{1}, b_{2} \})$ est ind\'ecomposable en consid\'erant la partition $%
p_{N_{m-1} \cup \{b_{2}\}}$. Il suffit de remarquer que $N_{m-1} \cup
\{b_{2}\} \longrightarrow m \longrightarrow b_{1}$ et que $b_{1} \not\in
[N_{m-1} \cup \{b_{2}\}]$. Il s'ensuit, d'apr\`es de corollaire 2.4, que $G
= G(N_{m}\cup\{b_{1}, b_{2}\})$. Ainsi $G-\{m, b_{1}\}$ est
ind\'ecomposable, contradiction car $b_{1}$ est un sommet isol\'e de $I(G)$.
{\hspace*{\fill}$\Box$\medskip}

\subsubsection{Les graphes (-1)-critiques $G$ en un sommet interne de $%
I^{\prime}(G)$}

Nous introduisons, pour $m \geq 2$ et $a\in\{1,\cdots, m-1\}$ , la classe $%
\mathcal{G}_{m}(a)$ des graphes $G$ d\'efinis sur $N_{m}$ tels que $%
N_{m}-\{0, 1\}$ et $N_{m-2}$ sont des intervalles respectifs de $G-0$ et de $%
G-m$; $N_{m}-\{1\}$ et $N_{m}-\{m- 1\}$ ne sont pas des intervalles de $G$
et pour tout $i\in\{1,\cdots, m-1\}-\{a\}$, $\{i-1, i+1\}$ est un intervalle
de $G-i$.

Notons d'abord les remarques suivantes.

\begin{remarque}
\'Etant donn\'e $G \in \mathcal{G}_{m}(a)$, la permutation $i
\longmapsto m-i$ est un isomorphisme de $G$ sur un graphe de $\mathcal{G}
_{m}(m-a)$.
\end{remarque}

\begin{remarque}
Soit $G\in \mathcal{G}_{m}(a)$. On a $G- \{i, i+1\} \simeq G - \{0, 1\}$ ou $%
G(N_{m-2})$ suivant que $0 \leq i \leq a-1$ ou que $a \leq i \leq m-1$
respectivement. Lorsque $m \geq 4$, $G(N_{m-2}) \in \mathcal{F}_{m-2}$ ou $%
\mathcal{G}_{m-2}(a)$ suivant que $a \geq m-2$ ou que $a < m-2$
respectivement; l'application $i \mapsto m-i$ est un isomorphisme de $G -
\{0,1\}$ sur un graphe de $\mathcal{F}_{m-2}$ ou de $\mathcal{G} _{m-2}(m-a)$
suivant que $a \leq 2$ ou que $a > 2$ respectivement. De plus, si $G - \{0,
1\}$ (resp. $G(N_{m-2})$) est ind\'ecomposable, $G$ est ind\'ecomposable si
et seulement si $0 \not \in [\{2, \cdots, m\}]$ (resp. $m \not \in [N_{m-2}]$%
).
\end{remarque}

\begin{lemme}
Soit $G$  un graphe d'ordre $\geq 7$, (-1)-critique en $a \in \{1,
\cdots, m-1\}$ et tel que $I^{\prime}(G) = P_{m}$. Alors, $G(N_{m})\in
\mathcal{G}_{m}(a)$.
\end{lemme}

\noindent\emph{Preuve\/. } On pose $G^{\prime}= G(N_{m})$. D'apr\`es le
lemme 2.1, pour tout $i\in\{1, \cdots, m-1\}-\{a\}$, $\{i-1, i+1\}$ est
un intervalle de $G-i$ et donc de $G^{\prime}-i$. De plus, $N_{m} -\{0,1\}$
et $N_{m-2}$ sont des intervalles respectifs de $G-0$ et de $G-m$. En
particulier, il s'agit d'intervalles respectifs de $G^{\prime}-0$ et de $%
G^{\prime}- m$. Il s'ensuit que $(1,0) \not \equiv (1,2) $ et que $(m-1,m-2)
\not \equiv (m-1,m)$, de sorte que $N_{m}-\{1\}$ et $N_{m} - \{m-1\}$ ne
sont pas des intervalles de $G^{\prime}$. Ainsi, $G^{\prime}\in \mathcal{G}%
_{m}(a)$. {\hspace*{\fill}$\Box$\medskip}

Nous introduisons les classes $\Omega_{1}$, $\Omega_{2}$ et $\Omega_{3}$ de
graphes d\'efinies comme suit.

\begin{itemize}
\item $\Omega_{1}$ est la classe des graphes $G$ d'ordre $\geq 7$,
(-1)-critiques en $2k+1$, tels que $I^{\prime}(G) = P_{2n+1}$, o\`u $n \geq
1 $ et $k \in N_{n-1}$.

\item $\Omega_{2}$ est la classe des graphes $G$ d'ordre $\geq 7$,
(-1)-critiques en $2k+1$, tels que $I^{\prime}(G) = P_{2n}$, o\`u $n \geq 1$
et $k \in N_{n-1}$.

\item $\Omega_{3}$ est la classe des graphes $G$ d'ordre $\geq 7$,
(-1)-critiques en $2k$, tels que $I^{\prime}(G) = P_{2n}$, o\`u $n \geq 2$
et $k \in \{1, \cdots, n-1\}$.
\end{itemize}

Soit $G$ un graphe d'ordre $\geq 7$, (-1)-critique en $2k$ et tel que $%
I^{\prime}(G) = P_{2n+1}$, o\`u $n \geq 1$ et $k \in \{1, \cdots, n\}$. La
permutation $\sigma$ de $S(G)$ qui fixe les sommets de $S(G) - N_{2n+1}$ et
telle que $\sigma(i) = 2n+1-i$ pour $i \in N_{2n+1}$, est un isomorphisme de
$G$ sur un graphe de la classe $\Omega_{1}$. Nous en d\'eduisons la remarque
suivante.

\begin{remarque}
\`A des isomorphismes pr\`es, les graphes $G$ d'ordre $\geq 7$,
(-1)-critiques en $a$ et tels que $I^{\prime}(G)$ est un chemin dont $a$ est
un sommet interne, sont les graphes de la classe $\Omega_{1} \cup \Omega_{2}
\cup \Omega_{3}$.
\end{remarque}


\

Nous caract\'erisons la classe $\mathcal{G}_{2n+1}(2k+1)$ comme suit.

\begin{lemme}
Soit un graphe $G$ d\'efini sur $N_{2n+1}$, o\`u $n \geq 1$. Alors $G\in \mathcal{G}%
_{2n+1}(2k+1)$, o\`u $k\in N_{n-1}$, si et seulement si $(0, 1)\not\equiv
(2, 1) \not\equiv (2n,2n+ 1)$ et pour tous $i \leq j \in N_{n}$ on a: pour $%
i \leq k$, $(2i, 2j+1)\equiv (0, 1)$; pour $i \geq k+1$, $(2i, 2j+1)\equiv
(2n, 2n+1)$; pour $i< j$, $(2i+1, 2j)\equiv (1, 2)\equiv (2i+1, 2j+1)$, $%
(2i, 2j)\equiv (0, 2)$ si $j \leq k$ et $(2i,2j) \equiv (1,2)$ si $j\geq k+1$%
.
\end{lemme}

\noindent\emph{Preuve\/. } Soient $n\geq 1$, $k\in N_{n-1}$ et $G\in
\mathcal{G}_{2n+1}(2k+1)$. Consid\'erons $i\leq j \in N_{n}$. Comme pour
tout $l\in \{1,\cdots,2n\}-\{2k+1\}$, $\{l-1, l+1\}$ est un intervalle de $%
G-l$ alors, pour $i\leq k$, $(2i, 2j+1)\equiv (0, 2j+1 )\equiv (0, 1)$; pour
$i\geq k+1$ , $(2i, 2j+1)\equiv (2i, 2n+1 )\equiv (2n, 2n+1 )$ et pour $%
i<j\leq k$, $(2i, 2j)\equiv (0, 2j )\equiv (0, 2)$. Comme de plus $%
N_{2n+1}-\{0,1\}$ et $N_{2n+1}-\{2n, 2n+1\}$ sont des intervalles respectifs
de $G-0$ et de $G-\{2n+1\}$, alors, pour $i < j$, $(2i+1, 2j)\equiv (1, 2j
)\equiv (1, 2 )\equiv (1, 2j+1 )\equiv (2i+1, 2j+1 )$ et, pour $j \geq k+1$
et $i<j$, $(2i, 2j)\equiv (2i, 2n )\equiv (1, 2n )\equiv (1, 2 )$. En outre,
$(0, 1)\not\equiv (2, 1)$ et $(2n, 2n+1)\not\equiv (2n, 2n-1 )\equiv (2, 1 )$
car $N_{2n+1}-\{1\}$ et $N_{2n+1}-\{2n\}$ ne sont pas des intervalles de $G$.

R\'eciproquement, on observe que pour tout $i\in \{2,\cdots,2n+1\}$ (resp. $%
i \in N_{2n-1}$), on a $(1, i )\equiv (1, 2 )$ (resp. $(i, 2n )\equiv (1, 2 )
$), c'est-\`a-dire, $\{2,\cdots,2n+1\}$ et $N_{2n-1}$ sont des intervalles
respectifs de $G-0$ et de $G-\{2n+1\}$. Comme de plus $(2n, 2n+1 )\not\equiv
(2, 1 )\not\equiv (0, 1 )$, $N_{2n+1}-\{1\}$ et $N_{2n+1}-\{2n\}$ ne sont
pas des intervalles de $G$. Enfin, pour $i\in\{1,\cdots,2n\}-\{2k+1\}$ et
pour $x\in N_{2n+1}-\{i-1,i, i+1\}$, on v\'erifie que $(x, i-1 )\equiv (x,
i+1 )$ de sorte que $\{i-1,i+1\}$ est un intervalle de $G-i$.
{\hspace*{\fill}$\Box$\medskip}

Nous caract\'erisons maintenant les graphes ind\'ecomposables de $\mathcal{G}
_{2n+1}(2k+1)$.

\begin{lemme}
Soient $n\geq 1$, $k \in N_{n-1}$ et $G\in \mathcal{G}_{2n+1}(2k+1)$. Le
graphe $G$ est ind\'ecomposable si et seulement si $(0, 1)\not \equiv (1, 2)$%
.
\end{lemme}

\noindent\emph{Preuve\/. }
Soit un graphe $G\in \mathcal{G}_{2n+1}(2k+1)$. Supposons que $(0, 1) \equiv
(1, 2)$. D'apr\`es le lemme 2.17, $(0, 1) \not\equiv (2, 1)$. Ainsi, quitte
\`a remplacer $G$ par $G^{\star}$, on peut supposer que $0\longrightarrow 1$%
. On v\'erifie \`a l'aide du lemme 2.17, que $N_{2k+1}
\longrightarrow\{2k+2,\cdots,2n+1\}$, en particulier $G$ est
d\'ecomposable. Supposons maintenant que $(0, 1)\not \equiv (1, 2)$ et
montrons, par r\'ecurrence sur $n$, que pour tout $k \in N_{n-1}$, les
graphes de $\mathcal{G}_{2n+1}(2k+1)$ sont ind\'ecomposables. Avec le lemme
2.17, on v\'erifie que les graphes de $\mathcal{G}_{3}(1)$ sont
ind\'ecomposables. Soit $n > 1$ et soit $G\in\mathcal{G}_{2n+1}(2k+1)$, o\`u
$k \in N_{n-1}$. D'apr\`es les remarques 2.14 et 2.13, $G-\{2n, 2n+1\}\in
\mathcal{G}_{2n-1}(2k+1)$ lorsque $k\neq n-1$; $G-\{0, 1\}$ est isomorphe
\`a un graphe de $\mathcal{G}_{2n-1}(2n-3)$ lorsque $k = n-1$. Il s'ensuit,
par hypoth\`ese de r\'ecurrence, que $G-\{2n, 2n+1\}$ ou $G-\{0, 1\}$ est
ind\'ecomposable. De plus, $2n+1 \not \in [N_{2n-1}]$ et $0 \not \in [\{2,
\cdots, 2n+1\}]$. Il s'ensuit que $G$ est ind\'ecomposable d'apr\`es la
remarque 2.14.
{\hspace*{\fill}$\Box$\medskip}

Nous introduisons maintenant la classe $\mathcal{G}$ des graphes $G = (S,A)$
tels que $S = N_{2n+1}$ ou $N_{2n+1} \cup \{\alpha\}$, o\`u $n\geq 1$; $%
G(N_{2n+1}) \in \mathcal{G}_{2n+1}(2k+1)$ et $k \in N_{n-1}$; $(0,1) \equiv
(1,2)$ si et seulement si $S \neq N_{2n+1}$ et tels que:

\begin{itemize}
\item Si $S = N_{2n+1}$, alors $(2n, 2n+1) \not\equiv (0,1)$ lorsque $(0,2)
\equiv (1,2)$; $(2n, 2n+1) \not\equiv (1,2)$ lorsque $k=0$; $(0,1)
\not\equiv (0,2)$ lorsque $k=n-1$.

\item Si $S - N_{2n+1} = \{\alpha\}$, alors pour tout $i \in N_{n}$, $(2i+1,
\alpha) \equiv (1,2)$, $(2i, \alpha) \equiv (0, \alpha)$ si $i \leq k $, $%
(2i, \alpha) \equiv (1,0)$ si $i \geq k+1$; $(0, \alpha) \not\equiv (1,2) $;
$(0, \alpha) \equiv (\alpha,0)$ lorsque $(0,2) \equiv (1,2) \equiv (2n, 2n+1)
$.
\end{itemize}

\begin{lemme}
Les graphes d'ordre $\geq 7$ de la classe $\mathcal{G}$ sont des graphes de
la classe $\Omega_{1}$.
\end{lemme}

\noindent\emph{Preuve\/. } Soit $G=(S, A)$ un graphe d'ordre $\geq 7$ de la
classe $\mathcal{G}$, il existe $n \geq 2$, $k\in N_{n-1}$ tels que $S =
N_{2n+1}$ ou $N_{2n+1} \cup \{\alpha\}$ et $G(N_{2n+1}) \in \mathcal{G}%
_{2n+1}(2k+1) $. D'apr\`es la remarque 2.14, pour tout $i \in N_{2n}$, il
existe un isomorphisme $\sigma$ de $G(N_{2n+1})-\{i, i+1\} $ sur $%
G(N_{2n+1})-\{0, 1\}$ ou $G(N_{2n-1})$. Lorsque $S = N_{2n+1}\cup \{\alpha\}$, $\sigma$ se prolonge en un isomorphisme, fixant $\alpha$, de $G-\{i, i+1\}$
sur $G-\{0, 1\}$ ou $G-\{2n, 2n+1\}$. Pour tout $x \in S- \{1, 2n\}$, $G -\{x,0\}$ et $G- \{x,2n+1\}$ sont d\'ecomposables. De m\^eme si $\alpha\in S$, $G- \{\alpha, 2k+1\}$ est d\'ecomposable car $N_{2k}$ et $\{2k+2, \cdots,
2n+1\} $ en sont des intervalles. Il suffit de montrer que les graphes $G$, $%
G-\{0,1\}$ et $G- \{2n,2n+1\}$ sont ind\'ecomposables. En effet, dans ce
cas, par constuction de $\mathcal{G}$, tous les sommets de $G-\{2k+1\}$,
sont des sommets critiques de $G$. Il s'ensuit en utilisant ce qui
pr\'ec\`ede, que $I^{\prime}(G) = P_{2n+1}$ d'apr\`es le lemme 2.1. De plus,
$\{2k, 2k+2\}$ n'est pas un intervalle de $G-\{2k+1\}$. En effet, si $(0,
2)\equiv (1, 2)\not\equiv (0, 1) $, alors $S=N_{2n+1}$ et $(2k, 2n+1)\equiv
(0, 1)\not\equiv (2n, 2n+1)\equiv (2k+2, 2n+1)$. Supposons que $(0, 2)\equiv
(1, 2)\equiv (0, 1) $. Dans ce cas $S = N_{2n+1} \cup \{\alpha\}$. Si $(0,
2)\equiv (2n, 2n+1)$ (resp. $(0, 2)\not\equiv (2n, 2n+1)$), alors $(\alpha,
2k+2)\equiv(1, 2)\not\equiv(\alpha, 0)\equiv(\alpha, 2k)$ (resp. $(2k,
2n+1)\equiv (0, 1)\not\equiv (2n, 2n+1)\equiv (2k+2, 2n+1)$). Si $(0,
2)\not\equiv (1, 2) $, alors $k\neq 0$ et $(0, 2k+2)\equiv (1, 2)\not\equiv
(0, 2)\equiv (0, 2k)$. D'apr\`es le lemme 2.1, $2k+1$ est un sommet non
critique de $G$, et donc $G$ est (-1)-critique en $2k+1$.

Montrons pour finir, que $G$, $G-\{0,1\}$ et $G- \{2n,2n+1\}$ sont
ind\'ecomp\-osables. Supposons d'abord que $S= N_{2n+1}$. Comme $(0, 1)\not
\equiv (1, 2)$ alors, d'apr\`es le lemme 2.18, $G$ est ind\'ecomposable. Si $%
k \neq n-1$, d'apr\`es la remarque 2.14, $G - \{2n,2n+1\} \in \mathcal{G}
_{2n-1}(2k+1)$ et donc, $G - \{2n,2n+1\}$ est ind\'ecomposable d'apr\`es le
lemme 2.18. Si $k=n-1$, $G - \{2n,2n+1\} \in\mathcal{F}_{2n-1} $ et, comme $%
(0, 2)\not \equiv (0, 1)\not \equiv (1, 2)$, $G - \{2n,2n+1\}$ est
ind\'ecomposable d'apr\`es le lemme 2.10. L'application $\tau: l\longmapsto
l-2$ est un isomorphisme de $G -\{0, 1\}$ sur un graphe $H$ de $\mathcal{G}%
_{2n-1}(2k-1)$ ou de $\mathcal{F}_{2n-1}$, suivant que $k \geq 1$ ou que $k=0
$ respectivement. Si $k\geq 1$, alors $(2, 3)\not\equiv_{G} (3, 4)$,  et
donc $(0, 1)\not\equiv_{H} (1, 2)$. Si $k= 0$, alors $(2, 4)\not\equiv_{G}
(2, 3)\not\equiv_{G} (3, 4)$ et donc $(0, 2)\not\equiv_{H} (0,
1)\not\equiv_{H} (1, 2)$. D'apr\`es les lemmes 2.18 et 2.10, le graphe $H$,
et donc $G -\{0, 1\}$, est ind\'ecomposable. Supposons enfin que $S-
N_{2n+1} = \{\alpha\}$. Pour tout $i \in N_{n}$, on pose $X_{i}= N_{2i+1}
\cup \{\alpha\}$. On v\'erifie que $G(X_{0})$ est ind\'ecomposable. Soit $%
i\in N_{n-1}$. Supposons que $G(X_{i})$ est ind\'ecomposable. Si $i< k$,
alors $2i+3\in X_{i}(2i+1)$, $2i+2\not\in X_{i}(2i+1)$ et $(2i+2, 2i+1)\not
\equiv (2i+2, 2i+3)$. Si $i\geq k$, alors $2i+2\in [X_{i}]$, $2i+3\not\in
[X_{i}]$ et $(2i+2, 2i+1)\not \equiv (2i+2, 2i+3)$. Ainsi, $G(X_{i+1})$ est
ind\'ecomposable d'apr\`es le lemme 1.2. En particulier, $G-\{2n, 2n+1\}$ et
$G$ sont ind\'ecomposables. Montrons enfin que $G-\{0, 1\} $ est
ind\'ecomposable. Pour $k\geq 1$, l'application $\tau: l \longmapsto l-2$ se
prolonge en un isomorphisme, fixant $\alpha$, de $G - \{0,1\}$ sur un graphe
$G^{\prime}$. On v\'erifie que $G^{\prime}\in \mathcal{G}$ et on d\'eduit,
d'apr\`es ce qui pr\'ec\`ede, que $G^{\prime}$, et donc $G - \{0,1\}$ est
ind\'ecomposable. Supposons maintenant que $k=0$. Si $(2n, 2n+1) \equiv (1,
2)$, on v\'erifie que $G-\{0, 1\}$ est isomorphe au tournoi critique $%
V_{2n+1}$. Supposons alors que $(2n, 2n+1)\not \equiv (1, 2)$. Comme $(2,
4)\not \equiv (2, 3)\not \equiv (3, 4)$, d'apr\`es le lemme 2.10 et la
remarque 2.14, $G-\{0, 1, \alpha\}$ est ind\'ecomposable. On a $\alpha\in
Ext(Y)$, o\`u $Y= S-\{0, 1, \alpha\}$. En effet, $\alpha\not\in[Y]$ car $(2,
\alpha) \equiv (2, 1)\not \equiv (1, 2)\equiv (3, \alpha)$. De plus, pour
tout $i\in\{1,\cdots, n\}$, $(2i, 2n+1) \equiv (2n, 2n+1)\not \equiv (2,
1)\equiv (\alpha, 2n+1)$ et $(\alpha, 2)\equiv (1, 2)\not \equiv (2i+1, 2)$.
Il s'ensuit que pour tout $j\in\{2,\cdots, 2n+1\}$, $\alpha\not\in Y(j)$.
{\hspace*{\fill}$\Box$\medskip}

\begin{proposition}
\`A des isomorphismes pr\`es, les graphes de la classe $\Omega_{1}$ sont les
graphes d'ordre $\geq 7$ de $\mathcal{G}$.
\end{proposition}

\noindent\emph{Preuve\/. } Soit $G=(S, A)$ un graphe d'ordre $\geq 7$,
(-1)-critique en $2k+1$, tel que $I^{\prime}(G) = P_{2n+1}$, o\`u $n \geq 1$
et $k \in N_{n-1}$. D'apr\`es le lemme 2.15, $G(N_{2n+1})\in\mathcal{G}%
_{2n+1}(2k+1)$. D'apr\`es le lemme 2.18 et le corollaire 2.4, $%
S-N_{2n+1}\neq \emptyset$ si et seulement si $(0, 1) \equiv (1, 2)$. Supposons d'abord que $S=N_{2n+1}$. Si $(0, 2)\equiv (1, 2)$,
alors $(2n, 2n+1)\not \equiv (0, 1)$. Autrement, on v\'erifie \`a
l'aide du lemme 2.17, que $\{2k, 2k+2\}$ est un intervalle de
$G-\{2k+1\}$, contradiction. Si $k=0$ , alors $(0, 2)\equiv (1, 2)$.
On a $(2n, 2n+1) \not \equiv (1, 2)$ sinon $G -\{0, 1\} \simeq
\mathcal{O}_{2n-1}$, contradiction. Si $k= n-1$, $G-\{2n, 2n+1\}$
est, d'apr\`es la remarque 2.14, un graphe de $\mathcal{F}_{2n-1}$.
Comme ce graphe est ind\'ecomposable, alors $(0, 2)\not \equiv (0,
1)$ d'apr\`es le lemme 2.10. Supposons maintenant que
$S-N_{2n+1}\neq \emptyset$. \`A un isomorphisme pr\`es, $\alpha \in
S-N_{2n+1}$. D'apr\`es le lemme 2.1, pour tout $i \in
N_{n}$, $(2i+1, \alpha) \equiv (1,2)$, $(2i, \alpha) \equiv (0, \alpha)$ si $%
i \leq k$; $(2i, \alpha) \equiv (2,1)\equiv (1,0)$ si $i \geq k+1$. Par les
lemmes 2.17 et 2.15, $(2, 1) \not \equiv (1, 2)$. Ainsi, quitte \`a
remplacer $G$ par $G^{\star}$, on peut supposer que $0\longrightarrow 1$ et $%
1\longrightarrow 2$. Il existe $\gamma \in S-N_{2n+1}$ tel que $(0,
\gamma)\not \equiv (1, 2)$, sinon on a une contradiction en v\'erifiant que $%
(S-N_{2n+1})\cup\{2k+2,\cdots,2n+1 \}$ est un intervalle non trivial de $G$.
Si $(0,2) \equiv (1,2) \equiv (2n,2n+1)$, on v\'erifie, \`a l'aide du lemme
2.17, que $G(N_{2n+1}) =\mathcal{O}_{2n+1}$. Il s'ensuit que $(0, \gamma)
\equiv (\gamma, 0)$, autrement $G(N_{2n+1} \cup \{\gamma\}) \simeq V_{2n+3}$%
, contradiction d'apr\`es le corollaire 2.4. Ainsi $G(N_{2n+1} \cup
\{\gamma\})$ est isomorphe \`a un graphe de $\mathcal{G}$. Lorsque $n=1$, on
v\'erifie que $G(N_{2n+1} \cup \{\gamma\}) = G(\{0,1,2,3,\gamma\})$ est
ind\'ecomposable. Lorsque $n \geq 2$, $G(N_{2n+1} \cup \{\gamma\})$ est
ind\'ecomposable d'apr\`es le lemme 2.19. D'apr\`es le corollaire 2.4, $%
\gamma = \alpha$ et $S = N_{2n+1} \cup \{\gamma\}$.
{\hspace*{\fill}$\Box$\medskip}


\begin{lemme}
Soit un graphe $G$ d\'efini sur $N_{2n}$, o\`u $n \geq 1$. Alors $G\in \mathcal{G}%
_{2n}(2k+1)$, o\`u $k \in N_{n-1}$, si et seulement si $(0, 1)\not\equiv (2,
1)$ et pour tous $x < y \in N_{2n}$ on a: si $x$ et $y$ ne sont pas tous
deux pairs, alors $(x,y) \equiv (1,2)$; si $x$ et $y$ sont pairs, alors $%
(x,y) \equiv (0,2)$, $(2n-2,2n)$ ou $(2k,2k+2)$, suivant que $y \leq 2k$, $x
\geq 2k+2$ ou que $x \leq 2k < y$ respectivement.
\end{lemme}

\noindent\emph{Preuve\/. } Soient $n\geq 1$, $k \in N_{n-1}$ et supposons $G\in \mathcal{G}_{2n}(2k+1)$. On a $(0, 1)\not\equiv (2, 1)$ car $%
N_{2n}-\{1\}$ n'est pas un intervalle de $G$ et $N_{2n}-\{0,1\}$ est un
intervalle de $G-0$. Soient $x < y$ dans $N_{2n}$. Si $x$ et $y$ sont pairs,
comme pour tout $l\in \{1,\cdots,2n-1\}-\{2k+1\}$, $\{l-1, l+1\}$ est un
intervalle de $G-l$, alors $(x,y) \equiv (0,2)$, $(2n-2,2n)$ ou $(2k,2k+2)
$, suivant que $y \leq 2k$, $x \geq 2k+2$ ou que $x \leq 2k < y$
respectivement. Sinon, comme $N_{2n}-\{0, 1\}$ et $N_{2n-2}$ sont des
intervalles respectifs de $G-0$ et de $G - 2n$, alors $(x,y) \equiv (1,y)
\equiv (1,2)$ si $x$ est impair; $(x,y) \equiv (x,2n-1) \equiv (1,2n-1)
\equiv (1,2)$ si $y$ est impair.

R\'eciproquement, comme pour tout $i\in N_{2n}-\{0,1\}$ (resp. $i\in N_{2n-2}
$), $(1, i )\equiv (1, 2 )$ (resp. $(i, 2n-1 )\equiv (1, 2)$) alors $%
N_{2n}-\{0, 1\}$ et $N_{2n-2}$ sont des intervalles respectifs de $G-0$ et
de $G-2n$. Comme $(1, 2)\equiv (0, 1 )\not\equiv (2, 1 )$, alors quitte \`a
remplacer $G$ par $G^{\star}$, on peut supposer que $0\longrightarrow
1\longrightarrow 2$ et donc $N_{2n-2} \longrightarrow 2n-1\longrightarrow 2n$%
. En particulier, $N_{2n}-\{1\}$ et $N_{2n}-\{2n-1\}$ ne sont pas des
intervalles de $G$. Soit $i\in \{1,\cdots, 2n-1\}-\{2k+1\}$. Remarquons que
pour tout $x\in N_{2n}-\{i-1,i, i+1\}$, $(x, i-1 )\equiv (x, i+1 )$ , de
sorte que $\{i-1,i+1\}$ est un intervalle de $G-i$.
{\hspace*{\fill}$\Box$\medskip}

\begin{lemme}
Soient $n\geq 1$, $k \in N_{n-1}$ et $G\in \mathcal{G}_{2n}(2k+1)$. Le
graphe $G$ est ind\'ecomposable si et seulement si $(2k, 2k+2)\not \equiv
(1, 2)$.
\end{lemme}

\noindent\emph{Preuve\/. } Comme d'apr\`es le lemme 2.21, $(1,2) \not \equiv
(2,1)$, quitte \`{a} remplacer $G$ par $G^{\star}$, on peut supposer que $1
\longrightarrow 2$. Si $(2k, 2k+2) \equiv (1, 2)$, alors, encore par le
lemme 2.21, on a $N_{2k} \longrightarrow \{2k+1, \cdots, 2n\}$ et donc $G$
est d\'ecomposable. Supposons que $(2k, 2k+2) \not \equiv (1, 2)$. Pour $n=1$%
, $G \in \mathcal{G}_{2}(1)$ et $(0,2) \not \equiv (1,2) \not \equiv (1,0)
\not \equiv (2,0)$, donc $G$ est ind\'ecomposable. Soit $n > 1$. Par la
remarque 2.14 et par hypoth\`ese de r\'ecurrence, $G - \{2n-1, 2n\}$ ou $G -
\{0,1\}$ est ind\'ecomposable. De plus, $2n \not \in [N_{2n-2}]$ et $0 \not
\in [\{2, \cdots, 2n\}]$. Il s'ensuit que $G$ est ind\'ecomposable d'apr\`es
la remarque 2.14.
{\hspace*{\fill}$\Box$\medskip}

Nous introduisons la classe $\mathcal{G^{\prime}}$ des graphes $G \in
\mathcal{G}_{2n}(2k+1)$, o\`u $n \geq 1$ et $k \in N_{n-1}$, tels que $(2k,
2k+2) \not\equiv (1,2)$; $(2n-2,2n) \not\equiv (0,2)$ si $(2k,2k+2) \equiv
(0,2)$; $(2n-2, 2n) \not\equiv (1,2)$ si $k=0$ et $(0,2) \not\equiv (1,2)$
si $k = n-1$.

Nous compl\'etons la remarque 2.13 comme suit.

\begin{remarque}
\'Etant donn\'e un graphe $G \in \mathcal{G^{\prime}}$, l'application $\phi:
x \longmapsto 2n-x$ est un isomorphisme de $G$ sur un graphe de $\mathcal{%
G^{\prime}}$.
\end{remarque}

\begin{lemme}
Les graphes d'ordre $\geq 7$ de la classe $\mathcal{G^{\prime}}$ sont des
graphes de la classe $\Omega_{2}$.
\end{lemme}

\noindent\emph{Preuve\/. } Soit $G= (S, A)$ un graphe d'ordre $\geq 7$ de la
classe $\mathcal{G}^{\prime}$. Il existe $n \geq 3$ et $k \in N_{n-1}$, tels
que $G \in \mathcal{G}_{2n}(2k+1)$. Comme $(2k, 2k+2)\not\equiv (1, 2)$, d'apr\`es le lemme 2.22, $G$ est ind\'ecomposable. D'apr\`es la remarque
2.14, pour tout $i \in N_{n-1}$, $G-\{i,i+1\} \simeq G-\{0, 1\}$ ou $%
G-\{2n-1, 2n\}$. Pour tout $x \in S - \{1, 2n-1\}$, Les graphes $G-\{x, 0\}$
et $G-\{x, 2n\}$ sont d\'ecomposables. Il suffit de montrer que $G-\{0,1\}$
et $G-\{2n-1,2n\}$ sont ind\'ecomposables. En effet, dans ce cas, par
construction de $\mathcal{G}^{\prime}$, tous les sommets de $G - \{2k+1\}$
sont des sommets critiques de $G$. Il s'ensuit, en utilisant ce qui
pr\'ec\`ede, que $I(G) = P_{2n}$ d'apr\`es le lemme 2.1. De plus, $\{2k,
2k+2\}$ n'est pas un intervalle de $G-\{2k+1\}$. En effet, si $(2k,
2k+2)\not\equiv (0, 2) $, alors $(0, 2k)\equiv (0, 2)\not\equiv(2k,
2k+2)\equiv (0, 2k+2)$. Si $(2k, 2k+2)\equiv (0, 2)$, alors $(2n, 2k)\equiv
(2k+2, 2k) \equiv (2, 0)\not\equiv (2n, 2n-2)\equiv (2n, 2k+2) $. D'apr\`es
le lemme 2.1, $2k+1$ est un sommet non critique de $G$, et donc $G$ est
(-1)-critique en $2k+1$. Montrons pour finir, que $G-\{2n-1, 2n\}$ et $%
G-\{0,1\} $ sont ind\'ecomposables. D'apr\`es la remarque 2.14, $G-\{ 2n-1,
2n\}\in \mathcal{G}_{2n-2}(2k+1)$ ou $\mathcal{F}_{2n-2}$ suivant que $k <
n-1$ ou que $k=n-1$ respectivement. D'une part $(2k, 2k+2)\not\equiv (1, 2)$%
, d'autre part $(0,1) \not \equiv (0,2) \not \equiv (1,2)$ lorsque $k=n-1$.
Il s'ensuit, en utilisant les lemmes 2.22 et 2.10, que $G-\{2n-1, 2n\}$ est
ind\'ecomposable. D'apr\`es la remarque 2.23, il existe un isomorphisme $\phi
$ de $G$ sur un graphe $H$ de $\mathcal{G}^{\prime}$, avec $\phi(0) = 2n$ et
$\phi(1) = 2n-1$. La restriction de $\phi$ \`a $S - \{0,1\}$ est un
isomorphisme de $G-\{0,1\}$ sur $H - \{2n-1, 2n\}$. Le graphe $H - \{2n-1,
n\}$ \'etant ind\'ecomposable d'apr\`es ce qui pr\'ec\`ede, il en est de
m\^eme pour $G-\{0,1\}$.
{\hspace*{\fill}$\Box$\medskip}

\begin{proposition}
$\Omega_{2}$ est la classe des graphes d'ordre $\geq 7$ de la classe $%
\mathcal{G}^{\prime}$.
\end{proposition}

\noindent\emph{Preuve\/. } Soit $G=(S, A)$ un graphe d'ordre $\geq 7$,
(-1)-critiques en $2k+1$, tel que $I^{\prime}(G) = P_{2n}$, o\`u $n \geq 1$
et $k \in N_{n-1}$. D'apr\`es le lemme 2.15, $G(N_{2n})\in \mathcal{G}%
_{2n}(2k+1)$. On a $S = N_{2n}$, autrement, d'apr\`es le lemme 2.1, pour $%
\gamma \in S- N_{2n}$, $(1, 2)\equiv (1, \gamma ) \equiv (3, \gamma) \equiv
\cdots \equiv (2n-1, \gamma) \equiv (2n-1, 1) \equiv (2, 1)$. Contradiction
car, d'apr\`es le lemme 2.21, $(1,2) \equiv (0,1) \not \equiv (2,1)$.
D'apr\`es le lemme 2.22, $(2k, 2k+2)\not\equiv (1, 2)$. Si $(2k, 2k+2)\equiv
(0, 2)$, alors $(2n-2, 2n)\not\equiv (0, 2)$. Sinon on a une contradiction
en v\'erifiant, \`a l'aide du lemme 2.21, que $\{2k, 2k+2\}$ est un
intervalle du graphe ind\'ecomposable $G-\{2k+1\}$. Si $k=0$, alors $%
(2n-2,2n)\not\equiv (1, 2)$. Autrement, en utilisant le lemme 2.21, $%
\{2,\cdots,2n-1\}$ est un intervalle non trivial du graphe ind\'ecomposable $%
G-\{0, 1\}$, contradiction. Lorsque $k=n-1$, l'application $i \longmapsto
2n-i$ est un isomorphisme de $G$ sur un graphe $H$ de $\Omega_{2}$ et on a $%
H \in \mathcal{G}_{2n}(1)$. D'apr\`es ce qui pr\'ec\`ede, $(2n-2,2n) \not
\equiv_{H} (1,2)$, donc $(2,0) \not \equiv_{G} (2n-1, 2n-2) \equiv_{G} (2,1)$%
. {\hspace*{\fill}$\Box$\medskip}


\begin{lemme}
Soit un graphe $G$ d\'efini sur $N_{2n}$, o\`u $n \geq 2$. Alors $G \in \mathcal{G}%
_{2n}(2k)$, o\`u $k \in \{1, \cdots, n-1\}$, si et seulement si $(0,1) \not
\equiv (2,1)\not \equiv(2n-1,2n)$ et pour tous $x < y \in N_{2n}$, on a: $%
(x,y) \equiv (0,2)$ si $x$ et $y$ sont pairs; $(x,y) \equiv (0,1)$ si $x$
est pair, $y$ est impair et $y < 2k$; $(x,y) \equiv (2n-1,2n)$ si $x$ est
impair, $y$ est pair et $x > 2k$; $(x,y) \equiv (1,2)$ dans le reste des cas.
\end{lemme}

\noindent\emph{Preuve\/. }
Soient $n\geq 2$, $k\in \{1,\cdots,n-1\}$, $G\in \mathcal{G}_{2n}(2k)$ et $x
< y \in N_{2n}$. Comme pour tout $l\in \{1,\cdots,2n-1\}-\{2k\}$, $\{l-1,
l+1\}$ est un intervalle de $G-l$, alors $(x,y) \equiv (0,2)$ si $x$ et $y
$ sont pairs; $(x,y) \equiv (0,1)$ si $x$ est pair, $y$ est impair et $y < 2k
$; $(x,y) \equiv (2n-1,2n)$ si $x$ est impair, $y$ est pair et $x > 2k$.
Comme de plus $N_{2n}-\{0, 1\}$ et $N_{2n-2}$ sont des intervalles
respectifs de $G-0$ et de $G-2n$, on v\'erifie que dans le reste des cas $%
(x,y) \equiv (1,2)$. Puisque $N_{2n}-\{1\}$ et $N_{2n}-\{2n-1\}$ ne sont pas
des intervalles de $G$, alors $(0, 1)\not\equiv (2, 1 )\equiv (2n-1,
2n-2)\not\equiv (2n-1, 2n)$.

R\'eciproquement, comme pour tout $i\in N_{2n}-\{0,1\}$ (resp. $i\in N_{2n-2}
$), on a $(1, i )\equiv (1, 2 )$ (resp. $(i, 2n-1 )\equiv (1, 2 )$), alors $%
N_{2n}-\{0, 1\}$ et $N_{2n-2}$ sont des intervalles respectifs de $G-0$ et
de $G-2n$. Comme $(2n-1, 2n)\not\equiv (2, 1 )\not\equiv (0, 1 )$, alors $%
N_{2n}-\{1\}$ et $N_{2n}-\{2n-1\}$ ne sont pas des intervalles de $G$. Soit $%
j\in \{1,\cdots, 2n-1\}-\{2k\}$. Pour $x\in N_{2n}-\{j-1,j, j+1\}$, $(x, j-1
)\equiv (x, j+1 )$, et donc $\{j-1,j+1\}$ est un intervalle de $G-j$.
{\hspace*{\fill}$\Box$\medskip}

\begin{lemme}
Soient $n\geq 2$, $k\in\{1, \cdots, n-1\}$ et $G\in\mathcal{G}_{2n}(2k) $. Alors $G$ est ind\'ecomposable si et seulement si $(0,2) \not \equiv (1,2)$.
\end{lemme}

\noindent\emph{Preuve\/. }
Si $(0,2) \equiv (1,2)$, d'apr\`es le lemme 2.26, $N_{2k}$ est un intervalle
non trivial de $G$. Donc $G$ est d\'ecomposable. Supposons que $(0,2) \not
\equiv (1,2)$. Pour $n=2$, $G \in \mathcal{G}_{4}(2)$ et on v\'erifie avec
l'aide du lemme 2.26, que $G$ est ind\'ecomposable. Soit $n > 2$. Par la
remarque 2.14 et par hypoth\`ese de r\'ecurrence, $G - \{2n-1, 2n\}$ ou $G -
\{0,1\}$ est ind\'ecomposable. De plus, $2n \not \in [ N_{2n-2}]$ et $0 \not
\in [\{2, \cdots, 2n\}]$, donc $G$ est ind\'ecomposable d'apr\`es la
remarque 2.14.
{\hspace*{\fill}$\Box$\medskip}

Nous consid\'erons enfin la classe $\mathcal{G}^{\prime\prime}$ des graphes $%
G = (S,A)$ tels que $S = N_{2n}$, $N_{2n} \cup \{\alpha\}$ ou $N_{2n} \cup
\{\alpha, \beta\}$, o\`u $n\geq 2$; $G(N_{2n})\in \mathcal{G}_{2n}(2k)$ o\`u
$k \in \{1, \cdots, n-1\}$; $(0,2) \equiv (1,2)$ si et seulement si $S -
N_{2n} \neq \emptyset$ et tels que:

\begin{itemize}
\item Si $S = N_{2n}$, alors $(2n-1,2n)\not\equiv (1,2)$ lorsque $(0,1)
\equiv (1,2)$; $(0,2) \not\equiv (2n-1,2n)$ lorsque $k=1$; $(0,2) \not\equiv
(0,1)$ lorsque $k=n-1$.

\item Si $S-N_{2n} \neq \emptyset$, alors pour tous $x \in N_{2n}$ et $%
\gamma \in S-N_{2n}$, $(x,\gamma) \equiv (0,\gamma)$ si $x$ est pair, $%
(x,\gamma) \equiv (1,2)$ ou $(2,1)$ si $x$ est impair et suivant que $x < 2k$
ou que $x > 2k$ respectivement. De plus, si $S-N_{2n} = \{\alpha\}$, alors $%
(2,1) \not \equiv (0, \alpha) \not \equiv (1,2)$. Si $S-N_{2n} = \{\alpha,
\beta\}$, alors $(\beta, \alpha) \not \equiv (0, \alpha) \equiv (1,2) \not
\equiv (2,1) \equiv (0, \beta)$.
\end{itemize}

La remarque suivante compl\`ete la remarque 2.13.

\begin{remarque}
Soit $G= (S,A)$ un graphe d\'efini sur $N_{2n}$, sur $N_{2n} \cup \{\alpha\}$
ou sur $N_{2n} \cup \{\alpha, \beta\}$ o\`u $n \geq 2$. Si $G \in \mathcal{G}%
^{\prime\prime}$, alors la permutation $h$ de $S$ d\'efinie par $h(x) = 2n-x$
si $x \in N_{2n}$, $h(\alpha) = \alpha$ lorsque $S - N_{2n}= \{\alpha\}$, $%
h(\alpha) = \beta$ et $h(\beta) = \alpha$ lorsque $S - N_{2n} = \{\alpha,
\beta\}$, est un isomorphisme de $G$ sur un graphe de $\mathcal{G}
^{\prime\prime}$.
\end{remarque}

\begin{lemme}
\'Etant donn\'e un graphe $G$ de la classe $\mathcal{G}^{\prime\prime}$, il
existe $n \geq 2$ et $k \in \{1, \cdots, n-1\}$ tels que $G$ est
(-1)-critique en $2k$ et $I^{\prime}(G) = P_{2n}$.
\end{lemme}

\noindent\emph{Preuve\/. }
Soit $G=(S, A)$ un graphe de la classe $\mathcal{G}^{\prime\prime}$. Il
existe $n\geq 2$ et $k\in \{1,\cdots,n-1\}$ tels que $S-N_{2n} = \emptyset$,
$\{\alpha\}$ ou $\{\alpha, \beta\}$ et $G(N_{2n})\in\mathcal{G}_{2n}(2k)$.
Nous montrons d'abord que les graphes $G$, $G-\{0,1\}$ et $G- \{2n-1,2n\}$
sont ind\'ecomposables. Supposons d'abord que $S= N_{2n}$. Comme $(0,2) \not
\equiv (1,2)$, d'apr\`es le lemme 2.27, $G$ est ind\'ecomposable. Si $k \neq
n-1$, d'apr\`es la remarque 2.14, $G - \{2n-1,2n\} \in \mathcal{G}_{2n-2}(2k)
$ et donc, $G - \{2n-1,2n\}$ est ind\'ecomposable d'apr\`es le lemme 2.27.
Si $k=n-1$, $G - \{2n-1,2n\} \in \mathcal{F}_{2n-2}$, et comme $(0,1) \not
\equiv (0,2) \not \equiv (1, 2)$, alors $G- \{2n-1,2n\}$ est
ind\'ecomposable d'apr\`es le lemme 2.10. Supposons maintenant que $%
S-N_{2n}\neq\emptyset$. Pour tout $i \in \{1, \cdots, n\}$, on pose $X_{i} =
N_{2i} \cup (S-N_{2n})$. On va montrer par r\'ecurrence que $G(X_{i})$ est ind\'ecomposable. C'est vrai pour $i = 1$. Soit $i \in \{1, \cdots, n-1\}$, on suppose $G(X_{i})$ ind\'ecomposable. Si $i< k$, alors $2i+2\in X_{i}(2i)$, $2i+1\not\in
X_{i}(2i)$ et $(2i+1, 2i)\not \equiv (2i+1, 2i+2)$. Si $i\geq k$ , alors $%
2i+1\in [X_{i}]$, $2i+2\not\in [X_{i}]$ et $(2i+1, 2i)\not \equiv (2i+1,
2i+2)$. Ainsi, $G(X_{i+1})$ est ind\'ecomposable d'apr\`es le lemme 1.2. En
particulier, $G-\{2n-1, 2n\}$ et $G$ sont ind\'ecomposables. Enfin, nous
d\'eduisons, \`a l'aide de la remarque 2.28, que $G-\{0,1\}$ est
ind\'ecomposable. En effet, il existe un isomorphisme $h$ de $G$ sur un graphe $H$ de $\mathcal{G}^{\prime\prime}$, avec $h(0) = 2n$ et $h(1) = 2n-1$. La restriction de $h$ \`a $S - \{0,1\}$ est un isomorphisme de $G-\{0,1\}$
sur $H - \{2n-1, 2n\}$. Le graphe $H - \{2n-1, n\}$ \'etant ind\'ecomposable
d'apr\`es ce qui pr\'ec\`ede, il en est de m\^eme pour $G-\{0,1\}$.

Soit $i \in N_{2n-1}$. D'apr\`es la remarque 2.14, il existe un isomorphisme
$\sigma$ de $G(N_{2n})-\{i, i+1\}$ sur $G(N_{2n})-\{0, 1\}$ ou $%
G(N_{2n})-\{2n-1, 2n\}$. Lorsque $S -N_{2n} \neq \emptyset$, $\sigma$ se
prolonge en un isomorphisme, fixant chaque sommet de $S-N_{2n}$, de $G-\{i,
i+1\}$ sur $G-\{0, 1\}$ ou $G-\{2n-1, 2n\}$. Il s'ensuit que $(i,i+1)$ est
un arc de $I(G)$. De plus, pour $x \in S- \{1, 2n-1\}$, $G -\{x,0\}$ et $G-
\{x,2n\}$ sont d\'ecomposables. Lorsque $S -N_{2n} = \{\alpha\}$, $G-\alpha$
et $G- \{\alpha, 2k\}$ sont d\'ecomposables car dans ce cas, $N_{2k} \sim
\{2k+1, \cdots, 2n\}$. Lorsque $S -N_{2n} = \{\alpha, \beta\}$, $G-\alpha$, $%
G-\beta$, $G- \{\alpha, 2k\}$, $G- \{\beta, 2k\}$ et $G - \{\alpha, \beta\}$
sont d\'ecomposables car dans ce cas, $N_{2k } \sim \{2k+1, \cdots,
2n,\alpha\}$ et $(N_{2k} \cup \{\beta\}) \sim \{2k+1, \cdots, 2n\}$. Il
s'ensuit que tous les sommets de $S-\{2k\}$ sont des sommets critiques de $G
$ et que $I^{\prime}(G) = P_{2n}$ d'apr\`es le lemme 2.1. De plus, $\{2k-1,
2k+1\}$ n'est pas un intervalle de $G-2k$. En effet, si $(1,2) \equiv (2,1)$%
, alors $(0, 2k-1)\not\equiv (0,2k+1)$. Si $(1,2) \not \equiv (2,1)$ et $%
S-N_{2n} \neq \emptyset$, alors  $(\alpha, 2k-1)\not\equiv (\alpha, 2k+1)$.
Supposons enfin que $(1,2) \not \equiv (2,1)$ et  $S=N_{2n}$, alors $(0,
2k-1)\not\equiv (0,2k+1)$ ou $(2n, 2k-1)\not \equiv (2n,2k+1)$ suivant que $%
(0, 1)\not \equiv (1,2)$ ou que $(0, 1)\equiv (1,2)$ respectivement. Par le
lemme 2.1, $2k$ est un sommet non critique de $G$. Ainsi, $G$ est
(-1)-critique en $2k$.
{\hspace*{\fill}$\Box$\medskip}

\begin{proposition}
\`A des isomorphismes pr\`es, les graphes de la classe $\Omega_{3}$ sont les
graphes d'ordre $\geq 7$ de $\mathcal{G}^{\prime\prime}$.
\end{proposition}

\noindent\emph{Preuve\/. } Soit $G=(S, A)$ un graphe d'ordre $\geq 7$,
(-1)-critique en $2k$, tel que $I^{\prime}(G) = P_{2n}$, o\`u $n \geq 2$ et $%
k \in \{1, \cdots, n-1\}$. D'apr\`es le lemme 2.15, $G(N_{2n})\in\mathcal{G}%
_{2n}(2k) $. D'apr\`es le lemme 2.27 et le corollaire 2.4, $(0,2) \equiv
(1,2)$ si et seulement si $S-N_{2n} \neq\emptyset$. Supposons que $S=N_{2n}$%
. Si $(0,1)\equiv (1,2)$, alors $(2n-1,2n) \not\equiv (1,2)$. Autrement,
d'apr\`es le lemme 2.26, $\{2k-1, 2k+1\}$ est un intervalle non trivial du
graphe ind\'ecomposable $G-\{2k\}$, contradiction. Si $k=n-1$ (resp. $k=1$)
alors $(0,2) \not\equiv (0,1)$ (resp. $(0,2) \not\equiv (2n-1,2n)$) car
sinon $\{1,\cdots,2n-2\}$ (resp. $\{2,\cdots,2n-1\}$) est un intervalle non
trivial du graphe ind\'ecomposable $G-\{2n-1, 2n\}$ (resp. $G-\{0, 1\}$),
contradiction. Supposons maintenant que $S- N_{2n}\neq\emptyset$. Soit $x
\in N_{2n}$ et $\gamma \in S-N_{2n}$. D'apr\`es le lemme 2.1, $(x,\gamma)
\equiv (0,\gamma)$ si $x$ est pair, $(x,\gamma) \equiv (1,2)$ ou $(2,1)$ si $%
x$ est impair et suivant que $x < 2k$ ou que $x > 2k$ respectivement. S'il
existe $\mu \in S-N_{2n}$ tel que $(1,2) \not \equiv (0, \mu) \not\equiv
(2,1)$, alors $G(N_{2n} \cup \{\mu\})$ est isomorphe \`a un graphe de $%
\mathcal{G}^{\prime\prime}$ et, d'apr\`es le corollaire 2.4, $G =
G(N_{2n}\cup \{\mu\})$. Sinon, $(1,2) \not \equiv (2,1)$ et $S - N_{2n} =
E_{1}\cup E_{2}$, o\`u $E_{1}= \{\gamma\in S-N_{2n}: (0, \gamma)\equiv (1,
2)\}$ et $E_{2}= \{\gamma\in S-N_{2n}: (0, \gamma)\equiv (2, 1)\}$. Il
existe $\alpha_{1}\in E_{1} $ et $\alpha_{2} \in E_{2}$ tels que $%
(\alpha_{2}, \alpha_{1}) \not \equiv(1,2)$ sinon $(E_{2} \cup N_{2k}) \sim
(E_{1} \cup \{2k+1, \cdots, 2n\})$, contradiction car $G$ est
ind\'ecomposable. Ainsi, $G(N_{2n} \cup \{\alpha_{1}, \alpha_{2}\})$ est
isomorphe \`a un graphe de $\mathcal{G}^{\prime\prime}$ et donc $%
G=G(N_{2n}\cup \{\alpha_{1}, \alpha_{2}\})$ d'apr\`es le corollaire 2.4.
{\hspace*{\fill}$\Box$\medskip}

\subsection{Les graphes (-1)-critiques $G$ tels que $I^{\prime}(G)$ est un
arbre \'etoil\'e}

Soient un entier $k \geq 3 $, $p_{1},\cdots, p_{k}$, $k$ entiers $\geq 2$ et
$i \in\{1,\cdots,k\}$. On pose $i_{0} = 0$, $S_{i_{p_{i}}}=\{i_{0},%
\cdots,i_{p_{i}}\}$ et $S(p_{1}, \cdots, p_{k}) = \bigcup\limits_{l=1}^{k}
S_{l_{p_{l}}}$. On d\'esigne par $P_{i_{p_{i}}}$ le chemin d\'efini sur $%
S_{i_{p_{i}}}$ par $A(P_{i_{p_{i}}})=\{(i_{l}, i_{h}), \mid l-h\mid = 1\}$
et on note $\mathcal{A}(p_{1},\cdots, p_{k})$ l'arbre $0$-\'etoil\'e
d\'efini sur $S(p_{1}, \cdots, p_{k})$ et dont les branches sont $%
P_{1_{p_{1}}},\cdots,P_{k_{p_{k}}}$. Pour $i \neq j \in \{1, \cdots, k\}$,
on consid\`ere l'application $f_{i_{p_{i}}, j_{p_{j}}}$ d\'efinie sur $%
S_{i_{p_{i}}} \cup S_{j_{p_{j}}}$ par $f_{i_{p_{i}}, j_{p_{j}}}(i_{l}) =
p_{i} - l$ et $f_{i_{p_{i}}, j_{p_{j}}}(j_{h}) = p_{i}+ h$.

\begin{remarque}
Si $G$ est un graphe (-1)-critique en $0$, tel que $I^{\prime}(G)= \mathcal{A}%
(p_{1},\cdots,p_{k})$ alors l'application $f_{i_{p_{i}}, j_{p_{j}}}$ est un
isomorphisme de $G(S_{i_{p_{i}}} \cup S_{j_{p_{j}}})$ sur un graphe $%
G^{\prime}$ de la classe $\mathcal{G}_{p_{i}+p_{j}}(p_{i})$.
\end{remarque}

Conform\'ement \`a la proposition 2.5, nous distinguons les cas, suivant que
l'arbre \'etoil\'e $I^{\prime}(G)$ admet ou n'admet pas une branche de
longueur impaire. Nous construisons alors pour $k$ entiers non nuls $n_{1}$,
$n_{2},\cdots,n_{k}$ o\`u $k \geq 3$, deux classes de graphes $\mathcal{H}%
(2n_{1}+1,2n_{2},\cdots, 2n_{k})$ et $\mathcal{H}(2n_{1},\cdots, 2n_{k})$,
comme suit.

\begin{enumerate}
\item $\mathcal{H}(2n_{1}+1,2n_{2},\cdots, 2n_{k})$ est la classe des
graphes $G$ d\'efinis sur $S(2n_{1}+1,2n_{2}, \cdots, 2n_{k})$ et qui
v\'erifient les conditions suivantes.

\begin{itemize}
\item Chacun des graphes $G(\{0, 1_{1}\})$ et $G(\{i_{1}, i_{2}\})$, o\`u $%
i\in \{2,\cdots,k\}$, est non vide.

\item Pour tous $x\neq y\in S(2n_{1}+1,2n_{2},\cdots, 2n_{k})$, on a:

\begin{itemize}
\item Si $(x, y)\in A_{1} = \{(1_{2l+1}, 1_{2j+1}): 0\leq l < j\leq n_{1}\}$
, alors $(x,y)\equiv (1_{1}, 1_{3})$.

\item Si $(x,y)\in A_{i}= \{(i_{2l+1}, i_{2j}): 0\leq l < j\leq n_{i}\}$,
o\`u $i\in\{2,\cdots,k\}$, alors $(x,y)\equiv (i_{1}, i_{2})$.

\item Si $(x, y)\in E\cup F$, o\`u $E=\{(i_{2j}, 1_{2l+1}): \ 2\leq i\leq k,
\ 0\leq j\leq n_{i}, \ 0\leq l\leq n_{1}\}$ et $F =\{(1_{2j}, 1_{2l+1}):
0\leq j\leq l \leq n_{1}\}$ alors $(x,y)\equiv (0, 1_{1})$.

\item Si $\{(x, y), (y, x)\} \bigcap \ (( \bigcup\limits_{i=1}^{k} A_{i})
\bigcup E \bigcup F )=\emptyset$, alors $G(\{x, y\})$ est vide.
\end{itemize}
\end{itemize}

\item $\mathcal{H}(2n_{1},\cdots, 2n_{k})$ est la classe des graphes $G$
d\'efinis sur $S(2n_{1},\cdots, 2n_{k})$ ou $S(2n_{1},\cdots, 2n_{k}) \cup
\{\gamma\}$, o\`u $\gamma \not \in S(2n_{1},\cdots, 2n_{k})$, et qui
v\'erifient les conditions suivantes.

\begin{itemize}
\item Pour tout $i \in \{1, \cdots, k\}$, le graphe $G(\{i_{1}, i_{2}\})$
est non vide.

\item Pour tous $i_{p} \neq j_{q}$ dans $S(2n_{1}, \cdots, 2n_{2})$, on a: $%
i_{p} \longleftrightarrow j_{q}$ si $p$ et $q$ sont pairs et $\gamma \not
\in S(G)$ ; $(i_{p},j_{q}) \equiv (i_{1}, i_{2})$ si $i = j$, $p$ est
impair, $q$ est pair et $p < q$; $i_{p}--j_{q}$ dans le reste des cas.

\item Lorsque $\gamma \in S(G)$, le graphe $G(\{\gamma, 0\})$ est non vide
et pour tout $i_{p} \in S(G) - \{\gamma\}$, $(\gamma, i_{p}) \equiv (\gamma,
0)$ si $p$ est pair, et $\gamma -- i_{p}$ si $p$ est impair.
\end{itemize}
\end{enumerate}

Soit $G$ un graphe de la classe $\mathcal{H}(p_{1},p_{2},\cdots, p_{k})$,
o\`u $p_{1} = 2n_{1}$ ou $2n_{1}+1$, $p_{r} = 2n_{r}$ pour $r \in \{2,
\cdots k\}$, et soient $i \neq j$ dans $\{1, \cdots, k\}$. Notons les remarques
suivantes.

\begin{remarque}
Pour tout $q \in \{0, \cdots, p_{i}-1\}$, l'application $g_{i_{p_{i}}, q}$
d\'efinie sur $S(G)$ par $g_{i_{p_{i}}, q}(x) = i_{l-2}$ si $x = i_{l}$ avec
$l \geq q+2$, et $g_{i_{p_{i}}, q}(x) = x$ sinon, est un isomorphisme de $G
- \{i_{q}, i_{q+1}\}$ sur $G - \{i_{p_{i}-1}, i_{p_{i}}\}$.
\end{remarque}

\begin{remarque}
L'application $f_{i_{p_{i}}, j_{p_{j}}}$ est un isomorphisme de $%
G(S_{i_{p_{i}}} \cup S_{j_{p_{j}}})$ sur un graphe de la classe $\mathcal{G}%
_{p_{i}+ p_{j}}(p_{i})$. De plus, si $\gamma \in S(G)$ alors $%
G(S_{i_{p_{i}}} \cup S_{j_{p_{j}}} \cup \{\gamma\})$ est isomorphe \`a un
graphe de $\mathcal{G}^{\prime\prime}$, en particulier, il s'agit d'un
graphe (-1)-critique.
\end{remarque}

\begin{lemme}
\'Etant donn\'e un graphe $G$ de la classe $\mathcal{H} (2n_{1}+1,2n_{2},
\cdots, 2n_{k})$, $G$ est (-1)-critique en $0$ et $I(G)= \mathcal{A}%
(2n_{1}+1,2n_{2},\cdots, 2n_{k})$.
\end{lemme}

\noindent\emph{Preuve\/. } Posons $n= \mid S(2n_{1}+1,2n_{2}, \cdots,
2n_{k})\mid $, $p_{1}=2n_{1}+1$ et $p_{i}=2n_{i}$ pour $i\in\{2,\cdots,k\}$.
Notons que $n$ est pair et $n \geq 8$. Nous montrons par r\'ecurrence sur $n$,
que les graphes $G$ et $G-\{i_{p_{i}},i_{p_{i}-1}\}$, o\`u $%
i\in\{1,\cdots,k\}$, sont ind\'ecomposables. Supposons d'abord que $G$ est
d'ordre $n=8$, c'est-\`a-dire $G \in \mathcal{H}(3,2,2)$. L'application $f$
d\'efinie par $f(2_{2}) = 0$ , $f(2_{1}) = 1$, $f(0) = 2$, $f(3_{1}) = 3$, $%
f(3_{2}) = 4$ et $f(1_{1}) = \alpha$, est un isomorphisme de $%
G-\{1_{2},1_{3}\} $ sur un graphe $G^{\prime}$. On v\'erifie que $G^{\prime}$
est un graphe de la classe $\mathcal{G}^{\prime\prime}$, en particulier $%
G-\{1_{2},1_{3}\} $ est ind\'ecomposable. D'apr\`es la remarque 2.33 et le
lemme 2.18, les graphes $G-\{2_{1},2_{2}\} $ et $G-\{3_{1},3_{2}\} $ sont
ind\'ecomposables. De plus, $3_{1} \in[S(G-\{3_{1},3_{2}\}) ]$, $3_{2}
\not\in [S(G-\{3_{1},3_{2}\})]$ et $(3_{1},3_{2})\not\equiv (3_{1}, 0)$,
donc $G$ est ind\'ecomposable. Supposons maintenant que $G$ est d'ordre $n
\geq 10$. Si pour tout $t \in \{1, \cdots, k\}$, $p_{t} \leq 3$,
c'est-\`a-dire $G \in \mathcal{H}(3,2,2,\cdots,2)$, alors $k \geq 4$ et $G -
\{k_{p_{k}}, k_{p_{k}-1}\} \in \mathcal{H}(p_{1}, \cdots, p_{k-1})$. Sinon,
il existe $t \in \{1, \cdots, k\}$ tel que $p_{t} \geq 4$ et donc $G -
\{t_{p_{t}}, t_{p_{t}-1}\} \in \mathcal{H}(q_{1}, \cdots, q_{k})$, o\`u $%
q_{t} = p_{t}-2$ et pour tout $r \in \{1, \cdots, k\} - \{t\}$, $q_{r} =
p_{r}$. Il s'ensuit en appliquant l'hypoth\`ese de r\'ecurrence, qu'il
existe $l \in \{1, \cdots, k\}$ tel que $G - \{l_{p_{l}}, l_{p_{l}-1}\}$ est
ind\'ecomposable, et tel que pour tout $m \in \{1, \cdots, k\} - \{l\}$, $G
- \{l_{p_{l}}, l_{p_{l}-1}, m_{p_{m}}, m_{p_{m}-1}\}$ est ind\'ecomposable.
On v\'erifie maintenant que $G - \{m_{p_{m}}, m_{p_{m}-1}\}$ et $G$ sont
ind\'ecomposables en utilisant le lemme 1.2 avec les parties $X = S(G) -
\{l_{p_{l}}, l_{p_{l}-1}, m_{p_{m}}, m_{p_{m}-1}\}$ et $Y = S(G) -
\{m_{p_{m}}, m_{p_{m}-1}\}$. En effet, $l_{p_{l}-1} \in [X]$, $l_{p_{l}}
\not \in [X]$ et $(l_{p_{l}-1}, l_{p_{l}}) \not \equiv (l_{p_{l}-1}, 1_{1})$%
, donc $G - \{m_{p_{m}}, m_{p_{m}-1}\}$ est ind\'ecomposable. De plus, $%
m_{p_{m}-1} \in [Y]$, $m_{p_{m}} \not \in [Y]$ et $(m_{p_{m}-1}, m_{p_{m}})
\not \equiv (m_{p_{m}-1}, 1_{1})$, donc $G$ est ind\'ecomposable.

Soit $j \in \{1, \cdots, k\}$. On a $S(G) - \{j_{p_{j}}, j_{p_{j}-1}\}$ est
un intervalle non trivial de $G - j_{p_{j}}$. De plus, pour tout $l \in \{1,
\cdots, p_{j}-1\}$, $\{j_{l-1}, j_{l+1}\}$ est un intervalle de $G - j_{l}$.
Il s'ensuit que tous les sommets de $G - {0}$ sont des sommets critiques de $%
G$. Comme $G-\{j_{p_{j}},j_{p_{j}-1}\}$ est ind\'ecomposable, les arcs de $%
\mathcal{A}(p_{1}, p_{2}, \cdots, p_{k})$ sont des arcs de $I(G)$ d'apr\`es
la remarque 2.32. Nous concluons, en utilisant le lemmes 2.1 et la
proposition 2.5, que $G$ est (-1)-critique en $0$ et que $I(G) =\mathcal{A}%
(p_{1}, p_{2}, \cdots, p_{k})$.
{\hspace*{\fill}$\Box$\medskip}

\begin{proposition}
Les graphes $G$ d'ordre $\geq 7$, (-1)-critiques en $0$ et tels que $%
I^{\prime}(G)=\mathcal{A}(2n_{1}+1,2n_{2},\cdots, 2n_{k})$ sont, aux
compl\'ementaires pr\`es, les graphes de la classe $\mathcal{H}%
(2n_{1}+1,2n_{2},\cdots, 2n_{k})$.
\end{proposition}

\noindent\emph{Preuve\/. } Soit $G=(S,A)$ un graphe (-1)-critique en $0$ tel que $I^{\prime}(G)=
\mathcal{A}(2n_{1}+1,2n_{2},\cdots, 2n_{k})$. Montrons que $G$ ou $\overline{%
G}$ est un graphe de $\mathcal{H}(2n_{1}+1,2n_{2},\cdots, 2n_{k})$. Notons
que pour $i\neq j$ dans $\{2,\cdots,k\}$, on a $(1_{2n_{1}},
1_{2n_{1}-1})\equiv(1_{2n_{1}}, i_{2n_{i}-1})\equiv(j_{2n_{j}-1},
i_{2n_{i}-1})\equiv(j_{2n_{j}-1}, 1_{2n_{1}})\equiv(1_{2n_{1}-1},
1_{2n_{1}}) $.

 Ainsi, le graphe $G(\{1_{2n_{1}}, 1_{2n_{1}-1}\})$ est
complet ou vide. Quitte \`a remplacer $G $ par $\overline{G}$, on peut
supposer que $1_{2n_{1}}--1_{2n_{1}-1}$. Nous v\'erifions \`a l'aide de la
remarque 2.31 et des lemmes 2.17 et 2.26, que $G(S(2n_{1}+1,2n_{2}\cdots,
2n_{k}))$ est un graphe de la classe $\mathcal{H} (2n_{1}+1,2n_{2}\cdots,
2n_{k})$. Ce graphe \'etant ind\'ecomposable d'apr\`es le lemme 2.32, \\$G =
G(S(2n_{1}+1,2n_{2},\cdots, 2n_{k}))$ d'apr\`es le corollaire 2.4.
{\hspace*{\fill}$\Box$\medskip}

\begin{lemme}
\'Etant donn\'e un graphe $G$ de la classe $\mathcal{H} (2n_{1},\cdots,
2n_{k})$, $G$ est (-1)-critique en $0$ et $I^{\prime}(G)= \mathcal{A}%
(2n_{1},\cdots, 2n_{k})$.
\end{lemme}

\noindent\emph{Preuve\/. } Posons $n= \mid S(G) \mid $ et $p_{i}=2n_{i}$
pour $i\in\{1,\cdots,k\}$. Notons que $n \geq 7$. Nous montrons par
r\'ecurrence sur $n$, que les graphes $G$ et $G-\{i_{p_{i}},i_{p_{i}-1}\}$,
o\`u $i\in\{1,\cdots,k\}$, sont ind\'ecomposables. Supposons d'abord que $G$
est d'ordre $n=7$, c'est-\`a-dire $G \in \mathcal{H}(2,2,2)$ et $S(G) =
S(2,2,2)$. D'apr\`es la remarque 2.33 et le lemme 2.27, les graphes $%
G-\{1_{1},1_{2}\}$, $G-\{2_{1},2_{2}\}$ et $G-\{3_{1},3_{2}\}$ sont
ind\'ecomposables. De plus, $3_{1}\in[S(2,2,2)-\{3_{1},3_{2}\}]$ , $3_{2}\not%
\in[ S(2,2,2)-\{3_{1},3_{2}\}]$ et $(3_{1},3_{2})\not\equiv (3_{1}, 0)$,
donc $G$ est ind\'ecomposable. Supposons maintenant que $G$ est d'ordre $%
n\geq 8$. Soit $i\in\{1,\cdots,k\}$. Si $p_{i} = 2$ et $k = 3$ alors,
d'apr\`es la remarque 2.33, $G - \{i_{p_{i}}, i_{p_{i}-1}\}$ est
ind\'ecomposable. Sinon, ou bien $p_{i} = 2$ et $k \geq 4$, ou bien $p_{i}
\geq 4$. Dans le premier cas, $G-\{i_{p_{i}},i_{p_{i}-1}\}$ est isomorphe
\`a un graphe de la classe $\mathcal{H}(q_{1}, \cdots, q_{k-1})$, o\`u pour
tout $t \in \{1, \cdots, k-1\}$, $q_{t} = p_{t}$ (resp. $p_{t+1}$) si $t < i$
(resp. si $t \geq i$). Dans le deuxi\`eme cas, $G \in \mathcal{H}(r_{1},
\cdots, r_{k})$, o\`u $r_{i} = p_{i} - 2$ et pour tout $t \in \{1, \cdots,
k\} - \{i\}$, $r_{t} = p_{t}$. Il s'ensuit dans les deux cas, que $%
G-\{i_{p_{i}},i_{p_{i}-1}\}$ est ind\'ecomposable par hypoth\`ese de
r\'ecurrence. Par construction de la classe $\mathcal{H} (2n_{1},\cdots,
2n_{k})$, $i_{p_{i}-1} \in [X]$, $i_{p_{i}} \not \in [X]$, o\`u $X = S(G) -
\{i_{p_{i}},i_{p_{i}-1}\}$, et $(i_{p_{i}-1}, i_{p_{i}}) \not \equiv
(i_{p_{i}-1}, 0)$, donc $G$ est ind\'ecomposable.

Soit $j \in \{1, \cdots, k\}$. Encore par construction de la classe $%
\mathcal{H} (2n_{1},\cdots, 2n_{k})$, pour tout $l\in\{1,\cdots,p_{j}-1\}$, $%
\{j_{l-1},j_{l+1}\}$ est un intervalle de $G-j_{l}$ et $S(G)
-\{j_{p_{j}-1},j_{p_{j}}\}$ est un intervalle non trivial de $G-j_{p_{j}}$.
De plus, si $\gamma \in S(G)$, alors $(S_{j_{p_{j}}}) -- (S(G-\gamma) -
S_{j_{p_{j}}})$, et donc $S_{j_{p_{j}}}$ est un intervalle non trivial de $G
- \gamma$. Il s'ensuit que tous les sommets de $G - 0$ sont des sommets
critiques de $G$. Comme $G-\{j_{p_{j}},j_{p_{j}-1}\}$ est ind\'ecomposable,
les arcs de $\mathcal{A}(p_{1}, p_{2}, \cdots, p_{k})$ sont des arcs de $I(G)
$ d'apr\`es la remarque 2.32. Il s'ensuit d'apr\`es le lemmes 2.1, que $G$
est (-1)-critique en $0$. De plus, lorsque $\gamma \in S(G)$, $S(G) -
\{\gamma, j_{p_{j-1}}, j_{p_{j}}\}$ est un intervalle non trivial de $G -
\{\gamma, j_{p_{j}}\}$, de sorte que d'apr\`es la proposition 2.5, $\gamma$
est un sommet isol\'e de $I(G)$. Encore par la la proposition 2.5, $%
I^{\prime}(G) =\mathcal{A}(p_{1}, p_{2}, \cdots, p_{k})$.
{\hspace*{\fill}$\Box$\medskip}

\begin{proposition}
Les graphes $G$ d'ordre $\geq 7$, (-1)-critiques en $0$ et tels que $%
I^{\prime}(G)=\mathcal{A}(2n_{1},2n_{2},\cdots, 2n_{k})$ sont, aux
compl\'ementaires pr\`es, les graphes de la classe $\mathcal{H}%
(2n_{1},2n_{2},\cdots, 2n_{k})$.
\end{proposition}

\noindent\emph{Preuve\/. }
Soit $G=(S,A)$ un graphe (-1)-critique en $0$ tel que $I^{\prime}(G)=
\mathcal{A}(2n_{1},\cdots, 2n_{k})$. Montrons que $G$ ou $\overline{G}$ est
un graphe de $\mathcal{H}(2n_{1},\cdots, 2n_{k})$. D'apr\`es la
remarque 2.31 et le lemme 2.26, pour $i \in \{2,\cdots,k\}$, $(i_{2n_{i}},
i_{2n_{i}-2}) \equiv (1_{2n_{1}-2}, 1_{2n_{1}})$ et $(i_{2n_{i}-1},
i_{2n_{i}-2}) \equiv (1_{2n_{1}-2}, 1_{2n_{1}-1})$. D'autre part, pour $i \neq j$ dans $\{2,\cdots,k\}$, $(1_{2n_{1}}, 1_{2n_{1}-2}) \equiv (1_{2n_{1}},
i_{2n_{i}}) \equiv (j_{2n_{j}},i_{2n_{i}}) \equiv (j_{2n_{j}},1_{2n_{1}})
\equiv (1_{2n_{1}-2}, 1_{2n_{1}})$ et $(1_{2n_{1}-1}, 1_{2n_{1}-2}) \equiv
(1_{2n_{1}-1}, i_{2n_{i}-1}) \equiv (j_{2n_{j}-1},i_{2n_{i}-1}) \equiv
(j_{2n_{j}-1},1_{2n_{1}-1}) \equiv (1_{2n_{1}-2}, 1_{2n_{1}-1})$. Les
graphes $G(\{1_{2n_{1}-2}, 1_{2n_{1}}\})$ et\\ $G(\{1_{2n_{1}-2},
1_{2n_{1}-1}\})$ sont alors vides ou complets. Quitte \`a remplacer $G$ par $%
\overline{G}$, on peut supposer qu'ou bien $1_{2n_{1}-2} \longleftrightarrow
1_{2n_{1}}$ et $1_{2n_{1}-2}--1_{2n_{1}-1}$, ou bien $1_{2n_{1}-2}--
\{1_{2n_{1}-1}, 1_{2n_{1}}\}$. Supposons d'abord que $1_{2n_{1}-2}
\longleftrightarrow 1_{2n_{1}}$ et $1_{2n_{1}-2}--1_{2n_{1}-1}$. Nous
v\'erifions en utilisant la remarque 2.31 et le lemme 2.26, que $%
G(S(2n_{1},\cdots, 2n_{k}))$ est un graphe de la classe $\mathcal{H}%
(2n_{1},\cdots, 2n_{k})$. Ce graphe \'etant ind\'ecomposable d'apr\`es le
lemme 2.36, $G = G(S(2n_{1},\cdots, 2n_{k}))$ d'apr\`es le corollaire 2.4.
Supposons maintenant que $1_{2n_{1}-2}-- \{1_{2n_{1}-1}, 1_{2n_{1}}\}$. On a
$0 - - (S(2n_{1},\cdots, 2n_{k})-\{0\})$. Comme $G$ est ind\'ecomposable il
existe $\gamma \in S(G)- S(2n_{1}, \cdots, 2n_{k})$ tel que le graphe $%
G(\{0, \gamma\})$ est non vide. En utilisant encore la remarque 2.31 et le
lemme 2.26, on v\'erifie que \\$G(S(2n_{1},\cdots, 2n_{k})\cup\{\gamma\})$ est
un graphe de la classe $\mathcal{H}(2n_{1},\cdots, 2n_{k})$. Ce graphe
\'etant ind\'ecomposable d'apr\`es le lemme 2.36, $G = G(S(2n_{1},\cdots,
2n_{k})\cup\{\gamma\})$ d'apr\`es le corollaire 2.4.
{\hspace*{\fill}$\Box$\medskip }

En conclusion nous obtenons le th\'eor\`eme suivant:

\begin{theoreme}
Les graphes $G$ d'ordre $\geq 7$ et $(-1)$-critiques sont, \`a isomorphisme pr\`es, les graphes $H_{2n+1}$, $\overline{H_{2n+1}}$, $R_{2n+1}$, $\overline{R_{2n+1}}$ o\`u $n \geq 3$; les graphes d'ordre $\geq 7$ de la classe $\mathcal{F} \cup \mathcal{G} \cup \mathcal{G'} \cup \mathcal{G''}$; les graphes de la classe $\mathcal{H}(2n_{1}+1, 2n_{2}, \cdots, 2n_{k}) \cup \mathcal{H}(2n_{1}, 2n_{2}, \cdots, 2n_{k})$ et leurs compl\'ementaires.

\end{theoreme}


Concernant les graphes $(-1)$-critiques d'ordre $\leq 6$, d'apr\`es le lemme 2.11, la classe $\mathcal{F}$ nous donne une famille de ces graphes. Remarquons alors que le lemme 2.3, outil important dans notre classification des graphes $(-1)$-critiques, ne s'\'etend pas \`a leur cas. Par exemple, le graphe $Q_5 = (\{0, 1, 2, \alpha, \beta\}, \{(0,1),(1,0),(0,2),(2,0),(0,\beta),(\beta,0),(2,\beta),(\beta,2),(\alpha,\beta),$\\$(\beta,\alpha)\})$ est un graphe $(-1)$-critique de la classe $\mathcal{F}$ dont le graphe \\d'ind\'ecomposabilit\'e est vide.

\end{document}